\documentclass[a4,12pt,epsfig]{elsarticle}
\usepackage{graphicx,epstopdf}

\usepackage{epsfig,epsf,epstopdf}
\usepackage{relsize}
\usepackage[pdf]{pstricks}
\usepackage{float}
\usepackage{mathptmx}
\usepackage{amsmath}
\usepackage{amssymb}
\usepackage{amsfonts}
\usepackage{mathtools}
\usepackage{euscript}
\usepackage{textcomp}
\usepackage{tikz, pgfplots}
\usetikzlibrary{shapes,arrows}
\usepackage{fancyhdr}
\usetikzlibrary{shadows,calc}
\pgfmathdeclarefunction{gauss}{2}{%
  \pgfmathparse{1/(#2*sqrt(2*pi))*exp(-((x-#1)^2)/(2*#2^2))}%
}

\newcommand\qedsymbol{\hbox{\rlap{$\sqcap$}$\sqcup$}}
\newcommand\smartqed{\renewcommand\qed{\relax\ifmmode\qedsymbol\else
  {\unskip\nobreak\hfil\penalty50\hskip1em\null\nobreak\hfil\qedsymbol
  \parfillskip=\z@\finalhyphendemerits=0\endgraf}\fi}}

\newcommand{\bea}{\begin{eqnarray}}
\newcommand{\eea}{\end{eqnarray}}
\newcommand{\nbea}{\begin{eqnarray*}}
\newcommand{\neea}{\end{eqnarray*}}

\newtheorem{theorem}{Theorem}[section]

\newtheorem{example}{Example}[section]

\begin{document}

\title{Reconstruction of Bivariate Piecewise Smooth Functions using Rational Approximation}
\author{ Akansha\footnote{E-mail: akansha.agrawal@manipal.edu}  \footnote{Work performed when author was at Indian Institute of Technology Bombay, Mumbai - 400076, India, E-Mail: akansha@math.iitb.ac.in} 
	\\
{ Department of Mathematics}\\{Manipal Institute of Technology Manipal}\\{Manipal Academy of Higher Education - 576104, India.} 
}
\date{}
\maketitle{}


\begin{center}{\bf Abstract}\end{center}
In this article, we develop a non-linear approximation technique to approximate piecewise smooth functions. We extend the concept of approximating piecewise smooth univariate functions using rational approximation into two-dimensional space. This article introduces the novel approach of piecewise Maehly-based Pad\`e-Chebyshev approximation. Initially, we develop a method called PiPC to approximate univariate piecewise smooth functions, then extend it to two-dimensional space, resulting in a bivariate piecewise Padé-Chebyshev approximation (Pi2DPC) for approximating piecewise smooth functions in two dimensions. We investigate the effectiveness of these techniques in minimizing the Gibbs phenomenon during the approximation of piecewise smooth functions. Specifically, we focus on approximating a special type of non-smooth functions with singularities along the axes. To the best of our knowledge, no previous work has addressed the approximation of non-smooth bivariate functions. Finally, we support our methods with numerical results to validate their ability to significantly reduce the Gibbs phenomenon.


\section{Introduction} 
Polynomials are classical tools to approximate a function. Runge \cite{run-85a} proved an exponential convergence of polynomial approximants of an analytic function $f$ in maximum norm on an approximation domain. Moreover, polynomial interpolation at $n$ Chebyshev points (roots of an $n$ degree Chebyshev polynomial $T_n(t)$ of first kind) converge to an analytic function $f$ at a rate $O(\rho^{-n})$, as $n\rightarrow \infty$, where $\rho>1$ is a constant (for details see \cite{riv-74a,boy-01a}). However, if $f$ is not analytic (or derivatives of $f$ has discontinuities), the exponential convergence of polynomial interpolations is not possible \cite{ber-12a}. Bernstein \cite{ber-14a} also proved that the global convergence of polynomial approximants of $|x|$ on $[-1,1]$ in maximum norm remains $O(1)$. The classical convergence rate deteriorates due to the presence of well known oscillatory behavior of the global polynomial approximants of functions which involves singularities\footnote{a point at which function is not defined or it is not differentiable, we refer to as a singularity.}. 

Polynomials serve as classical tools for function approximation. Runge \cite{run-85a} demonstrated exponential convergence of polynomial approximations of an analytic function $f$ in the maximum norm over an approximation domain. Additionally, polynomial interpolation at $n$ Chebyshev points (roots of an $n$-degree Chebyshev polynomial $T_n(t)$ of the first kind) converges to an analytic function $f$ at a rate of $O(\rho^{-n})$ as $n\rightarrow \infty$, where $\rho>1$ is a constant (for further details, refer to \cite{riv-74a,boy-01a}). However, when $f$ is not analytic or has discontinuous derivatives, achieving exponential convergence of polynomial interpolations is not feasible \cite{ber-12a}. Bernstein \cite{ber-14a} further demonstrated that the global convergence of polynomial approximations of $|x|$ on $[-1,1]$ in the maximum norm remains $O(1)$. The classical convergence rate diminishes due to the well-known oscillatory behavior of global polynomial approximations of functions involving singularities, defined as points where a function is either undefined or not differentiable.

In practice one often encounters the problem of approximating a function which is not smooth but piecewise smooth (function which involves singularities). In applied mathematics, physics, and engineering, the physical quantity of interest like, temperature, pressure, density, stress, strain etc. given by a scalar field $f$, needs to be evaluated over a certain region of space. Such an $f$, obtained as the solution of a certain PDE (governed by physical consideration of the problem), generally turns out to be piecewise smooth (for example, nonlinear hyperbolic PDEs). Approximating a smooth function is a well understood \cite{nur_89a,phi_03a,sho_73a,dev_lor-93a} subject of constructive approximation theory.  However, while approximating a piecewise smooth function, one has to face many difficulties. The main challenge in approximating piecewise smooth functions using polynomial interpolants, splines, truncated series etc., is the appearance of well-known \textit{Gibbs phenomenon} in the vicinity of a singularity. A discussion on approximating analytic but non-periodic functions defined on $[-1,1]$ using truncated Fourier series is discussed in \cite{lan-16a}, where the oscillatory behavior (Gibbs oscillations) of the approximants can be seen near $\pm 1$. Since the typical global approximation methods are suffered from Gibbs phenomenon while approximating a piecewise smooth function, there are methods which construct a piecewise polynomial approximation to recover such functions. A recursive numerical method is proposed by Trefthen \textit{et. at.} \cite{pac_pla_tre-10a} which automatically detects the location of singularities (of all orders) present in a bounded piecewise smooth function $f$ and construct a piecewise Chebyshev interpolation which recovers $f$ accurately. Eckhoff \cite{eck-93a} approximates the locations and the magnitude of a jump discontinuity of $f$ by solving a system of algebraic equations and reconstructs the function using step functions.  It can be seen in recent developments \cite{dmi_yos-12a,bat-15a} that such techniques works well in resolving the Gibbs phenomenon.

There are methods based on rational approximations to approximate functions with limited regularity. In fact, Newman \cite{new-64a} proved that a sequence of rational approximants of $|x|$ converges with a root exponential rate and the convergence rate can not be improved further. Furthermore, Pad\'e approximants which are rational approximations of a sequence of truncated power series, offer a faster rate of convergence (than the corresponding series approximants), rendering themselves a more suitable choice for approximating analytic functions (see Baker \cite{bak-geo-gra-pet_96a,bak-geo_65a}). Litvinov \cite{lit-03a} observe the property of error auto correction in Pad\'e approximations. A truncated Fourier series expansion of a periodic, even or odd, piecewise smooth function is used by Geer \cite{gee_ban-97a} to construct a Pad\'e approximant, which has been further developed by Min et al. \cite{min_kab_don-07a} for general functions. Pad\'e-Legendre interpolation method is proposed by Hesthaven et al. \cite{hes_kab_lur-06a} to approximate a piecewise smooth function. Kaber and Maday \cite{kab_mad-05a} obtained results on the convergence of a sequence of Pad\'e-Chebyshev approximations to the sign function.  Driscoll and Fornberg \cite{dri_for-01a} (also see Tampos \textit{et al.} \cite{tam_lop_hes-12a}) developed a singular Pad\'e approximation technique using finite series expansion coefficients or function values which approximates a piecewise smooth function much faster. These nonlinear methods minimize the Gibbs phenomenon significantly. However, these approximation techniques are limited to the univariate setting.

A more challenging problem is approximating a multivariate piecewise smooth function which involves jumps across some curve(s) in a bounded domain. There are several high accuracy algorithms to approximate univariate piecewise smooth functions (as mentioned previously) which have been used to improve the convergence rate of the approximants. However, generalization of these methods to more than one variable is not straightforward. Moreover, similar to the univariate case, the linear approximants (for example polynomial approximation, series partial sums, splines etc.) for multivariate functions with discontinuities exhibits sustained Gibbs oscillations in the vicinity of curve(s) of singularities. There is not much progress has been made in resolving the problem of Gibbs phenomenon in the case of multivariate piecewise smooth functions. Levin in \cite{lev-20a}  reconstructs non-periodic, smooth and non-smooth univariate functions as well as piecewise smooth bivariate functions using their Fourier series coefficients. His main idea involves identifying the location and the structure (for example in the case of 2-D functions) of the singularities, using B-spline functions and subsequently, approximating the smooth segments with higher order linear approximants. A different approximation method which does not require to approximate the locations of the singularities present in the function, using Christoffel-Darboux kernels is introduced by Marx et al. \cite{mar_pau_wei_hen_las-21a} to approximate multivariate piecewise smooth functions.

We develop a local approximation technique based on Pad\'e-Chebyshev approximation, most suitable for approximating univariate and multivariate piecewise smooth functions. The idea of Pad\'e-Chebyshev approximation based on series expansion other than \cite{dri_for-01a,tam_lop_hes-12a}, proposed by Maehly (see \cite{mae_60a}) and a more general case by Cheney (see \cite{che_66a}), is used in developing the approximation technique. We first propose a piecewise Pad\'e-Chebyshev (refered to as PiPC) algorithm for approximating univariate non smooth functions and then extend the work of Maehly to bivariate functions, in order to afford notional economy and since further generalization to higher dimensions is routine. We develop a global bivariate Pad\'e-Chebyshev approximation technique using the definition introduced by Lutterodt \cite{lut_74a}, wherein a specific choice of equations forming a linear system ascertains that the resulting coefficient matrix is the sum of block Toeplitz and block Hankel matrices; structurally similar to the univariate case. Having developed a global approximation method for approximating bivariate functions, we proceed towards a piecewise implementation of this method (henceforth referred to as Pi2DPC) to approximate bivariate non-smooth functions. This algorithm affords lesser computational complexity and the key advantage that it shares with PiPC is that, it suppresses the Gibbs phenomenon significantly without any prior knowledge of the type and location of the singularity present in the original function. We further supply numerical justification to validate the effectiveness of our algorithms in approximating piecewise smooth functions.

The outline of this article consisting of seven sections is as follows : In section \ref{sec:chebyshevseries} we briefly discuss the univariate Chebyshev series approximation. The basic construction of the global Pad\'e-Chebyshev method, suggested by Maehly, is recalled in section \ref{sec:MPCreconstruction} and subsequently, the PiPC algorithm is proposed. In Section \ref{sec:1Dnumericalresults}, we compare the numerical results obtained from PiPC with those from global Pad\'e-Chebyshev. We then propose a piecewise bivariate Chebyshev series approximation (henceforth referred to as Pi2DC) in section \ref{sec:2Dchebyshevapprox} and extend the global Pad\'e-Chebyshev approximation, defined in section \ref{sec:MPCreconstruction}, to a two dimensional space in section \ref{sec:2DMPCapprox}. Section \ref{sec:Pi2DPC} is devoted to the development of Pi2DPC algorithm. The article concludes by considering a few test examples to study the performance of the Pi2DPC method and making a comparison of the results obtained from Pi2DC and Pi2DPC with those from their global counterparts, in section \ref{sec:2Dnumericalresults}.

\section{Univariate Chebyshev series expansion}\label{sec:chebyshevseries}
For a given function $f$ which is Lipschitz continuous on $[-1,1]$, the Chebyshev series representation of $f$ is given by
\begin{equation}\label{eq:CheExp}
f(t)= \sum_{j=1}^{\infty}{\vphantom{\sum}}' c_jT_j(t)=:\mathsf{C}_\infty[f](t),
\end{equation}
where $T_j(t)$ denotes the Chebyshev polynomials of first kind of degree $j$, the prime in the summation indicates that the first term is halved and $c_j$ for $j=0,1,2,\ldots$ are the Chebyshev coefficients given by (see \cite{nur_89a,mas-han-03a})
\begin{equation}\label{eq:Checoff}
c_j=\dfrac{2}{\pi}\int_{-1}^{1}\dfrac{f(t)T_j(t)}{\sqrt{1-t^2}}dt.
\end{equation}
Computation of the integral in \eqref{eq:Checoff} for an arbitrary function is not always feasible hence, we use the Gauss-Chebyshev quadrature formula (see \cite{mas-han-03a,riv-74a})
\begin{equation}\label{eq:approxcoeff}
c_{j,n}:=\frac{2}{n}\sum_{l=1}^{n}f(t_l)T_k(t_l), \hspace{1cm} j=0,1,2,\ldots,
\end{equation}
where $t_l,l=1,2,\ldots,n$ are the Chebyshev points and are given by
\begin{equation}\label{eq:ChePts}
t_l	= \cos\left(\dfrac{\pi}{n}\left(l-\dfrac{1}{2}\right)\right), \hspace{1cm} l=1, 2, \ldots,n. 
\end{equation}  
The Chebyshev series expansion for a function defined on an arbitrary interval $[a,b]$ can be obtain using the affine map
\begin{equation}\label{eq:bijection}
\mathsf{G}(t) = a+(b-a)\frac{(t+1)}{2}, \hspace{1cm} t\in [-1,1].
\end{equation}
\section{Univariate Pad\'e-Chebyshev Reconstruction}\label{sec:MPCreconstruction}
The Pad\'e-Chebyshev approximation method based on Chebyshev series expansion of a function (other than \cite{dri_for-01a,tam_lop_hes-12a,aka_bas-19a}) is proposed by Maehly \cite{mae_60a} which we recall in this section. Let $f$ be a function defined on $[-1,1]$ such that it can be represented by the Chebsyhev series given in \eqref{eq:CheExp}. For the given integers $n_p \geq n_q \geq 1$, a rational function 
$$\mathsf{R}_{n_p,n_q}(t) := \dfrac{P_{n_p}(t)}{Q_{n_q}(t)} = \dfrac{{\displaystyle\sum_{i=0}^{n_p}}p_iT_i(t)}{{\displaystyle\sum_{j=0}^{n_q}}q_jT_j(t)}$$
with numerator polynomial $P_{n_p}(t)$ of degree  $\leq n_p$  and denominator polynomial $Q_{n_q}(t)$ of degree $\leq n_q$ with $Q_{n_q} \neq 0$, satisfying
\begin{equation}\label{eq:MaehlyPC}
Q_{n_q}(t)\mathsf{C}_\infty[f](t)-P_{n_p}(t)= \mathit{O}(T_{n_p+n_q+1}(t)), \hspace*{.5cm} t\rightarrow 0,
\end{equation}
is called a Pad\'e-Chebyshev approximant to the function $f$ of order $(n_p,n_q)$. Using the identity 
$$T_iT_j = \dfrac{1}{2}\left(T_{i+j}+T_{|i-j|}\right),$$
on the left hand side of \eqref{eq:MaehlyPC} and comparing the coefficients on both sides, we obtain the following homogeneous
\begin{equation}\label{1Dhomo}
\sum_{j=0}^{n_q}\left(c_{k-j}+c_{k+j}\right)q_j = 0, ~~~~~~~ k = n_p+1, n_p+2, \ldots, n_p+n_q+1,
\end{equation}
and the inhomogeneous system
\begin{align}\label{1Dinhomo}
\begin{split}
p_0 & = \dfrac{1}{2}\sum_{j=0}^{n_q}{\vphantom{\sum}}''c_jq_j\\
p_k &= \dfrac{1}{2}\sum_{j=0}^{n_q}(c_{k-j}+c_{k+j}+c_{j-k})q_j, ~~~~~~~ k = 1,2,\ldots,n_p,
\end{split}
\end{align}
where $c_j=0$ for $j<0$ and the double prime ($''$) indicates that the first term in the summation is doubled. Let us write the homogeneous system \eqref{1Dhomo} in the matrix form as
\begin{equation}\label{1Dhomomatrix}
\left(\begin{bmatrix}
c_{n_p+1} & c_{n_p}  & \cdots & c_{n_p-n_q+1}\\
c_{n_p+2} & c_{n_p+1} & \cdots & c_{n_p-n_q+2}\\
\vdots \\
c_{n_p+n_q} & c_{n_p+n_q-1} &  \cdots & c_{n_p}
\end{bmatrix}+\begin{bmatrix}
c_{n_p+1} & c_{n_p+2}  & \cdots & c_{n_p+n_q+1}\\
c_{n_p+2} & c_{n_p+3} & \cdots & c_{n_p+n_q+2}\\
\vdots \\
c_{n_p+n_q} & c_{n_p+n_q+1} &  \cdots & c_{n_p+2n_q}
\end{bmatrix}\right)\begin{bmatrix}
q_0\\q_1\\ \vdots \\ q_{n_q}\end{bmatrix} = \begin{bmatrix}
0\\0\\ \vdots \\ 0
\end{bmatrix}.
\end{equation}
The coefficient matrix in the above linear homogeneous system is the sum of a Toeplitz and a Hankel matrix of size $n_q\times (n_q+1)$ and therefore a non-trivial solution of \eqref{1Dhomomatrix} always exists.

We write the above system as
\begin{equation}\label{eq:1Dhomoshortmatrix}
(T_{n_p,n_q}^{Q}+H_{n_p,n_q}^Q)\textbf{q} = \textbf{0},
\end{equation}
where, for $r=n_p+1,n_p+2,\ldots,n_p+n_q, s=1,2,\ldots,n_q+1,$
$$T_{n_p,n_q}^Q = (t_{r,s}),~~t_{r,s} = c_{r-s+1}$$
is the Toeplitz matrix,
$$H_{n_p,n_q}^Q = (h_{r,s}),~~h_{r,s} = c_{r+s-1},$$ is the Hankel matrix 
and $\textbf{q}=(q_0,q_1,\ldots,q_{n_q})^T$.

Once the coefficient vector for the denominator polynomial is calculated, the numerator polynomial coefficients can be computed by the following matrix vector multiplication
\begin{multline}\label{1Dinhomomatrix}
\begin{bmatrix}
p_0\\p_1\\p_2\\ \vdots \\ p_{n_p}
\end{bmatrix} = \dfrac{1}{2}\left(\begin{bmatrix}
c_0 \\
c_1 & c_0\\
\vdots \\
c_{n_q} & c_{n_q-1}  & \cdots & c_0\\
\vdots \\
c_{n_p} & c_{n_p-1}  & \cdots & c_{n_p-n_q}
\end{bmatrix}+\right.\\
\left.\begin{bmatrix}
c_0 & c_1  & \cdots & c_{n_q}\\
c_1 & c_2  & \cdots & c_{n_q+1}\\
\vdots \\
c_{n_q} & c_{n_q+1}  & \cdots & c_{2n_q}\\
\vdots \\
c_{n_p} & c_{n_p+1} & \cdots & c_{n_p+n_q}
\end{bmatrix}+\begin{bmatrix}
0 & 0  & 0 &\cdots & 0 \\
0 & c_0  & c_1 &\cdots & c_{n_q-1} \\
\vdots \\
0 & 0  & 0 &\cdots & c_0 \\
\vdots \\
0 & 0  & 0 & \cdots & 0
\end{bmatrix}\right)\begin{bmatrix}
q_0\\q_1\\q_2\\ \vdots \\ q_{n_q}\end{bmatrix}.
\end{multline}

We denote the resulting Maehly based Pad\'e-Chebyshev approximant by $\mathsf{R}_{n_p,n_q}^\mathsf{M}(t)$. Note that, if the matrices in \eqref{eq:1Dhomoshortmatrix} is of full rank, then a normalization $q_0 = 1$ leads to a unique (normalized) $(n_p,n_q)$ order Pad\'e-Chebyshev approximant of $f$. 

Note that, $\mathsf{R}^{\mathsf{M}}_{n_p,n_q}(t)$ can be computed for a given set of Chebyshev coefficients $\{c_0, c_1, \ldots, c_{n_p+2n_q}\}$. However, in our numerical experiments we use the approximated coefficients, $c_{j,n}, j=0,1,\ldots,n_p+2n_q$, obtained by using $n$ Chebyshev points in the Gauss-Chebyshev quadrature formula \eqref{eq:approxcoeff}. The resulting approximant is denoted by  $\mathsf{R}^{n,\mathsf{M}}_{n_p,n_q}(t)$ and the coefficient matrix in \eqref{eq:1Dhomoshortmatrix} by $T_{n_p,n_q}^{n,Q}+H_{n_p,n_q}^{n,Q} =: A_{n_p,n_q}^{n,Q}$. We use $Null(A_{n_p,n_q}^{n,Q})$ in MATLAB to calculate the denominator coefficient vector \textbf{q} and in the case of numerical rank deficiency one can choose any one basis vector of $Null(A_{n_p,n_q}^{n,Q})$ as a solution of \eqref{1Dhomomatrix}. We choose the last column of $Null(A_{n_p,n_q}^{n,Q})$ in our numerical simulations (for the sake of consistency).

\section{Bivariate Chebyshev Series Expansion}\label{sec:2Dchebyshevapprox}
For $(x,y)\in D := [-1,1]\times[-1,1]$ and given non-negative integers $i$ and $j$, a bivariate Chebyshev polynomial is defined as (see \cite{bas_73a}) 
$$T_{ij}(x,y)=T_i(x)T_j(y),$$
where $T_i(x)$ and $T_j(y)$ are the $i$ and $j$ degree Chebyshev polynomials of first kind. 
\begin{theorem}
(Mason \cite{mas-han-03a}) Let $f(x,y)$ be a continuous function of bounded variation with one of the partial derivative is bounded in $D$. Then $f$ can be represented by the Chebyshev series as 
\begin{equation}\label{eq:bivariChebexpansion}
f(x,y)=\sum_{i = 0}^{\infty}\sum_{j = 0}^{\infty}c_{i,j}T_i(x)T_j(y)=:\mathsf{C}_\infty[f](x,y).
\end{equation}
\end{theorem}
A bivariate polynomial (see \cite{bas_73a})
\begin{equation*}\label{eq:trunbivariateChebexp}
\mathsf{C}_{d_x,d_y}[f](x,y) := \sum_{i = 0}^{d_x}\sum_{j = 0}^{d_y}c_{i,j}T_i(x)T_j(y),
\end{equation*}
of degree $d_x+d_y+1$, 
where
\begin{equation}\label{eq:bivarChebcoeff}
c_{i,j} =\dfrac{\epsilon_{i,j}}{\pi^2} \int_{-1}^{1}\int_{-1}^{1}\dfrac{f(x,y)T_i(x)T_j(y)}{\sqrt{1-x^2}\sqrt{1-y^2}}dxdy,
\end{equation}
and
\begin{align}\label{epsilon}
\epsilon_{i,j} & = \begin{cases}
1, ~~~~~ i = j = 0,\\
2, ~~~~~ i \ne j =0 ~~\mbox{ or }~~ j\ne i = 0, \nonumber \\
4, ~~~~~ i \ne j \ne 0.
\end{cases}
\end{align}
is called truncated bivariate Chebyshev series of $f$. Computing the coefficients \eqref{eq:bivarChebcoeff} accurately for an arbitrary function is not always feasible. Therefore, to approximate the coefficients for a given bivariate function $f$, we use the following two-dimensional Gauss-Chebyshev quadrature rule 
\begin{equation}\label{eq:bivarChebcoeffapprx}
c_{i,j} \approx \dfrac{\epsilon_{i,j}}{n_xn_y}\sum_{l_x=1}^{n_x}\sum_{l_y=1}^{n_y}f(t_{l_x},t_{l_y})T_i(t_{l_x})T_j(t_{l_y}) =: c_{i,j,n_x,n_y},
\end{equation}
for numerical integration in \eqref{eq:bivarChebcoeff}, where $t_{l_x}, l_x=1,2\ldots,n_x$ and $t_{l_y}, l_y=1,2\ldots,n_y$ are the Chebyshev points in $x$-direction and $y$-direction, respectively on a tensor grid. These points are the roots of the Chebyshev polynomials $T_{n_x}$ and $T_{n_y}$ given in \eqref{eq:ChePts} for $n = n_x$ and $n_y$, respectively.

Note that, we can compute a bivariate Chebyshev series approximant for a function defined on an arbitrary rectangular domain $[a_x,b_y]\times[a_y,b_y]$ using the affine map we defined in \eqref{eq:bijection}.
\subsection{Construction of Piecewise bivariate Chebyshev Approximation}\label{sec:Pi2DC}
Let $f$ be a bounded function defined on a rectangular domain $[a_x,b_x]\times[a_y,b_y]$. For given two-tuples $\mathbf{n}=(n_x,n_y)\ge \mathbf{1}$, $\mathbf{N}=(N_x,N_y)\ge \mathbf{1}$, and $(N_x-1)$-tuple $\mathbf{d_x}=\left(d_x^0,d_x^1,\ldots,d_x^{N_x-1}\right)$ and $(N_y-1)$-tuple $\mathbf{d_y}=\left(d_y^0,d_y^1,\ldots,d_y^{N_y-1}\right)$, construct the piecewise bivariate Chebyshev approximation of $f$ as follows:
\begin{enumerate}
\item Discretize the domain into $N_x$ cells in $x$-direction and $N_y$ cells in $y$-direction, denoted by $I_{j_x} = [a_{j_x},b_{j_x}],j_x = 0,1, \ldots, N_x-1$ and $I_{j_y} = [a_{j_y},b_{j_y}],j_y = 0,1,\ldots, N_y-1$, respectively, and denote the partition as $P_{N_x} := \{a_{0_x} , a_{1_x} ,\ldots, a_{N_x}\}$ and $P_{N_y} := \{a_{0_y} , a_{1_y} ,\ldots, a_{N_y}\}$, as explained in step 1 of Section \ref{subsec:PiPCM}.
\item Generate Chebyshev points $\{t_{l_x},l_x=1,2,\ldots,n_x\}$ and $\{t_{l_y},l_y=1,2,\ldots,n_y\}$ in the reference interval $[-1,1]$ using \eqref{eq:ChePts} with $n=n_x$ and $n=n_y$, respectively.  
\item For $j_x=0,1, \ldots, N_x-1$ and $j_y=0,1, \ldots, N_y-1$, obtain the approximated Chebyshev coefficients $c_{i,j,n}^{j_x\times j_y}$, for $i = 0, 1,\ldots, d_x^{j_x}$ and $j = 0,1,\ldots, d_y^{j_y}$, in each sub rectangle $I_{j_x}\times I_{j_y}$ using the formula \eqref{eq:bivarChebcoeffapprx} with the values of $f$ evaluated at $\mathsf{G}_{j_x}(t_{l_x})$ and $\mathsf{G}_{j_y}(t_{l_y})$. Here $\mathsf{G}_{j_x} : [-1,1] \rightarrow I_{j_x}$ and $\mathsf{G}_{j_y} : [-1,1] \rightarrow I_{j_y}$ are the maps defined in \eqref{eq:bijection}.
	\item The truncated Chebyshev series expansion of $f{|_{I_{j_x}\times I_{j_y}}}$ with $d_x^{j_x} d_y^{j_y} +1$ number of terms is given by    
	\begin{equation}\label{eq:2Dchexpcell}
	\mathsf{C}_{\mathbf{d_x,d_y}}^{n_x,n_y,I_{j_x}\times I_{j_y}}[f](x,y):=\sum_{i = 0}^{d_x^{j_x}}\sum_{j = 0}^{d_y^{j_y}}c_{i,j,n_x,n_y}T_i(x)T_j(y), 
	\end{equation}
	where  $(x,y) \in I_{j_x}\times I_{j_y}$, for all $j_x=0,1, \ldots, N_x-1$ and $j_y=0,1, \ldots, N_y-1$. 
\end{enumerate}
Define the piecewise $2$D Chebyshev approximation of $f$ in the rectangular domain $I_x\times I_y$, where $I_x= [a_x,b_x]$ and $I_y = [a_y,b_y]$, with respect to the given partition as
\begin{equation}
\mathsf{C}_{\mathbf{d_x,d_y}}^{n_x,n_y,N_x,N_y}[f](x,y) = 
\begin{cases}
\mathsf{C}_{d_x^0,d_y^0}^{n_x,n_y,I_{0_x}\times I_{0_y}}[f](x,y) &\quad\text{if } (x,y) \in [a_{0_x},b_{0_x})\times [a_{0_y},b_{0_y})\\
\mathsf{C}_{d_x^0,d_y^1}^{n_x,n_y,I_{0_x}\times I_{1_y}}[f](x,y) &\quad\text{if } (x,y) \in [a_{0_x},b_{0_x})\times [a_{1_y},b_{1_y})\\
\vdots \\
\mathsf{C}_{d_x^0,d_y^{N_y-1}}^{n_x,n_y,I_{0_x}\times I_{N_y-1}}[f](x,y) &\quad\text{if } (x,y) \in [a_{0_x},b_{0_x})\times [a_{N_y-1},b_{N_y-1})\\
\vdots \\
\mathsf{C}_{d_x^{N_x-1},d_y^0}^{n_x,n_y,I_{N_x-1}\times I_{0_y}}[f](x,y) &\quad\text{if } (x,y) \in [a_{N_x-1},b_{N_x-1})\times [a_{0_y},b_{0_y})\\
\mathsf{C}_{d_x^{N_x-1},d_y^1}^{n_x,n_y,I_{N_x-1}\times I_{1_y}}[f](x,y) &\quad\text{if } (x,y) \in [a_{N_x-1},b_{N_x-1})\times [a_{1_y},b_{1_y})\\
\vdots \\
\mathsf{C}_{d_x^{N_x-1},d_y^{N_y-1}}^{n_x,n_y,I_{N_x-1}\times I_{N_y-1}}[f](x,y) &\quad\text{if } (x,y) \in [a_{N_x-1},b_{N_x-1})\times [a_{N_y-1},b_{N_y-1})
\end{cases}
\end{equation} 
Hence $\mathsf{C}_{\mathbf{d_x,d_y}}^{n_x,n_y,N_x,N_y}[f](x,y)$ is the desired piecewise  bivariate Chebyshev approximation to the function $f$ on the domain $I_x\times I_y$ which we referred to as Pi2DC approximant of the function $f$.

\section{Bivariate Pad\'e-Chebyshev Approximation}\label{sec:2DMPCapprox} 
We extend the univariate Pad\'e-Chebyshev approximation method defined in Section \ref{sec:2DMPCapprox}, to two dimensional space. Our univariate case study on the Pad\'e-Chebyshev approximation suggests room for further improvement. Therefore, we implement a piecewise variant of the two-dimensional Pad\'e-Chebyshev approximation for more accurate results. In general (as observed in the univariate case), piecewise Pad\'e is a simple technique with less computational complexity to approximate a non-smooth function. Also, prior knowledge of the jump location, magnitude, or the order singularities present in the function to be approximated, is not required.

Consider a function $f$ defined on the unit square $D$ and represented by a double Chebyshev series \eqref{eq:bivariChebexpansion}. For given two-tuples $\mathbf{n_p} = (n_{p_x},n_{p_y})\ge\mathbf{n_q} = (n_{q_x},n_{q_y})\ge\mathbf{1}(=(1,1))$, where $n_{p_x},n_{q_x},n_{p_y},n_{q_y}$ are integers,  a rational function
\begin{equation}\label{rationalfunc}
\mathsf{R}_{\mathbf{n_p},\mathbf{n_q}}(x,y) = \dfrac{P_{\mathbf{n_p}}(x,y)}{Q_{\mathbf{n_q}}(x,y)},~~~~ (x,y) \in D,
\end{equation}
with numerator polynomial
\begin{equation}\label{numpoly}
P_{\mathbf{n_p}}(x,y) = \sum_{i = 0}^{n_{p_x}}\sum_{j = 0}^{n_{p_y}} p_{i,j}T_i(x)T_j(y),
\end{equation}
and denominator polynomial
\begin{equation}\label{denpoly}
Q_{\mathbf{n_q}}(x,y) = \sum_{r = 0}^{n_{q_x}}\sum_{s = 0}^{n_{q_y}} q_{r,s}T_r(x)T_s(y) \ne 0,
\end{equation}
of degrees $\le\mathbf{n_p}$ and $\le\mathbf{n_q}$, respectively, is called a bivariate Pad\'e-Chebyshev approximation of $f$, if
\begin{equation}\label{2Dpade}
\mathsf{C}_\infty[f](x,y)Q_{\mathbf{n_q}}(x,y) - P_{\mathbf{n_p}}(x,y) = O(T_{n_{p_x}+n_{q_x}+1}(x)T_{n_{p_y}+n_{q_y}+1}(y)), \hspace{.5cm} x\rightarrow 0, y\rightarrow 0.
\end{equation}
To compute an approximant $\mathsf{R}_{\mathbf{n_p},\mathbf{n_q}}(x,y)$, we need to determine the coefficients of the polynomials $P_{\mathbf{n_p}}(x,y)$ and $Q_{\mathbf{n_q}}(x,y)$, that is, we need to calculate $\tau := (n_{p_x}+1)(n_{p_y}+1)+(n_{q_x}+1)(n_{q_y}+1)$ unknowns. To compute these unknown parameters, we use the following definition (see \cite{lut_74a,lut_75a,you-chu_09a} for details) 
\begin{equation}\label{polycoeffconditions}
\mathsf{C}_\infty[f](x,y)Q_{\mathbf{n_q}}(x,y) - P_{\mathbf{n_p}}(x,y) = \sum_{i,j = 0}^{\infty}\eta_{i,j}T_i(x)T_j(y),
\end{equation}
where $\eta_{i,j}=0$ for
$$(i,j) \in P:=\{(i,j):0\leq i\leq n_{p_x},0\leq j\leq n_{p_y}\},$$
and
$$(i,j) \in J\setminus\{n_{p_x}+n_{q_x}+1,n_{p_y}+n_{q_y}+1\},$$
when
$$J:=\{(i,j) : n_{p_x}+1 \leq i \leq n_{p_x}+n_{q_x}+1, n_{p_y}+1 \leq j \leq n_{p_y}+n_{q_y}+1\}.$$
Observe that, the relation \eqref{polycoeffconditions} relates the polynomials $P_{\mathbf{n_p}}(x,y)$ and $Q_{\mathbf{n_q}}(x,y)$ and gives $\tau-1$ linear equations in $\tau$ unknowns which can be split into a homogeneous 
\begin{equation}\label{2Dhomosystem}
\sum_{r = 0}^{n_{q_x}}\sum_{s = 0}^{n_{q_y}}(c_{i-r,j-s}+c_{i-r,j+s}+c_{i+r,j-s}+c_{i+r,j+s})q_{r,s} = 0,
\end{equation}
for $(i,j)\in J\setminus\{(n_{p_x}+n_{q_x}+1,n_{p_y}+n_{q_y}+1)\}$, and an inhomogeneous part
\begin{align}\label{2Dinhomosystem}
p_{0,0} &= \dfrac{1}{4}\left(c_{0,0}q_{0,0}+\sum_{s = 0}^{n_{q_y}}c_{0,s}q_{0,s}+\sum_{r = 0}^{n_{q_x}}c_{r,0}q_{r,0}+\sum_{r = 0}^{n_{q_x}}\sum_{s = 0}^{n_{q_y}}c_{r,s}q_{r,s}\right),\nonumber\\
p_{i,0}& = \dfrac{1}{4}\left(\sum_{r = 0}^{n_{q_x}}(c_{i-r,0}+c_{i+r,0}+c_{r-i,0})q_{r,0}+\sum_{r = 0}^{n_{q_x}}\sum_{s = 0}^{n_{q_y}}(c_{i-r,s}+c_{i+r,s}+c_{r-i,s})q_{r,s}\right),~~~~~~ i = 1,2,\ldots,n_{p_x},\nonumber\\
p_{0,j} &=\dfrac{1}{4}\left(\sum_{s = 0}^{n_{q_y}}(c_{0,j-s}+c_{0,j+s}+c_{0,s-j})q_{0,s}+ \sum_{r = 0}^{n_{q_x}}\sum_{s = 0}^{n_{q_y}}(c_{r,j-s}+c_{r,j+s}+c_{r,s-j})q_{r,s}\right),~~~~~~ j = 1,2,\ldots,n_{p_y},\nonumber\\
p_{i,j}&=\dfrac{1}{4} \sum_{r = 0}^{n_{q_x}}\sum_{s = 0}^{n_{q_y}}(c_{i-r,j-s}+c_{i-r,j+s}+c_{i-r,s-j}+c_{i+r,j-s}+c_{r-i,j-s}+c_{i+r,j+s}+c_{i+r,s-j}+c_{r-i,j+s}+c_{r-i,s-j})q_{r,s},
\end{align}
for $(i,j)\in P$, where $c_{i,j}=0,i,j<0$. A further expansion of indexes in \eqref{2Dhomosystem} and \eqref{2Dinhomosystem} leads to block matrix systems. 

Let us denote the denominator and the numerator coefficient vectors as
\begin{equation*}
\textbf{Q}_{vec}:= \begin{bmatrix}
q_0^{n_{q_y}+1}\\
q_1^{n_{q_y}+1}\\
q_2^{n_{q_y}+1}\\
\vdots \\
q_{n_{q_x}}^{n_{q_y}+1}
\end{bmatrix} ~~ \mbox{ and } ~~ \textbf{P}_{vec}:= \begin{bmatrix}
p_0^{n_{p_y}+1}\\
p_1^{n_{p_y}+1}\\
p_2^{n_{p_y}+1}\\
\vdots \\
p_{n_{p_x}}^{n_{p_y}+1}
\end{bmatrix},
\end{equation*}
respectively where ${q}_{s}^{n_{q_y}+1}:= \left(q_{0,s},q_{1,s},\ldots,q_{n_{q_x},s}\right)^T$ and ${p}_{v}^{n_{p_y}+1}:= \left(p_{0,v},p_{1,v},\ldots
p_{n_{p_x},v}\right)^T$ for $s = 0,1,2,\ldots, n_{q_y}$ and $v = 0,1,2,\ldots, n_{p_y}$

On substituting $i = n_{p_x}+1,n_{p_x}+2,\ldots, n_{p_x}+n_{q_x}+1$ and $j=n_{p_y}+1,n_{p_y}+2,\ldots, n_{p_y}+n_{q_y}+1$ in \eqref{2Dhomosystem}, we get the following homogeneous linear system 
\begin{multline*}
\left(\begin{bmatrix}
T_{n_{p_x}+1} & T_{n_{p_x}} &  \cdots & T_{n_{p_x}-n_{q_x}+1}\\
T_{n_{p_x}+2} & T_{n_{p_x}+1} &  \cdots & T_{n_{p_x}-n_{q_x}+2}\\
T_{n_{p_x}+3} & T_{n_{p_x}+2} &  \cdots & T_{n_{p_x}-n_{q_x}+3}\\
\vdots & & & &\\
T_{n_{p_x}+n_{q_x}+1} & T_{n_{p_x}+n_{q_x}} & \cdots & T_{n_{p_x}+1}
\end{bmatrix}+\begin{bmatrix}
H_{n_{p_x}+1} & H_{n_{p_x}} &  \cdots & H_{n_{p_x}-n_{q_x}+1}\\
H_{n_{p_x}+2} & H_{n_{p_x}+1} &  \cdots & H_{n_{p_x}-n_{q_x}+2}\\
H_{n_{p_x}+3} & H_{n_{p_x}+2} &  \cdots & H_{n_{p_x}-n_{q_x}+3}\\
\vdots & & & &\\
H_{n_{p_x}+n_{q_x}+1} & H_{n_{p_x}+n_{q_x}} & \cdots & H_{n_{p_x}+1}
\end{bmatrix}+\right.\\
\left.\begin{bmatrix}
T_{n_{p_x}+1} & T_{n_{p_x}+2} & \cdots & T_{n_{p_x}+n_{q_x}+1}\\
T_{n_{p_x}+2} & T_{n_{p_x}+3} & \cdots & T_{n_{p_x}+n_{q_x}+2}\\
T_{n_{p_x}+3} & T_{n_{p_x}+4}  & \cdots & T_{n_{p_x}+n_{q_x}+3}\\
\vdots & & &\\
T_{n_{p_x}+n_{q_x}+1} & T_{n_{p_x}+n_{q_x}+2} & \cdots & T_{n_{p_x}+2n_{q_x}+1}
\end{bmatrix}+\begin{bmatrix}
H_{n_{p_x}+1} & H_{n_{p_x}+2} & \cdots & H_{n_{p_x}+n_{q_x}+1}\\
H_{n_{p_x}+2} & H_{n_{p_x}+3} & \cdots & H_{n_{p_x}+n_{q_x}+2}\\
H_{n_{p_x}+3} & H_{n_{p_x}+4}  & \cdots & H_{n_{p_x}+n_{q_x}+3}\\
\vdots & & &\\
H_{n_{p_x}+n_{q_x}+1} & H_{n_{p_x}+n_{q_x}+2} & \cdots & H_{n_{p_x}+2n_{q_x}+1}
\end{bmatrix}\right)\mathbf{Q}_{vec}=\mathbf{0},
\end{multline*}
where each element is a matrix of size $(n_{q_y}+1)\times(n_{q_y}+1)$ defined as 
$$T_m(u,v):=c_{m,u-v} \hspace{1cm} \mbox{ and } \hspace{1cm} H_m(u,v):=c_{m,u+v},$$
for $u=n_{p_y}+1,n_{p_y}+2,\ldots, n_{p_y}+n_{q_y}+1, v=0,1,2,\ldots,n_{q_y}$ and $m = n_{p_x}-n_{q_x}+1, \ldots, n_{p_x}+1,n_{p_x}+2,\ldots,n_{p_x}+2n_{q_x}+1$. 

We can write the above defined block matrix linear system in the following shorthand form
\begin{equation}\label{2DComputingQ}
(T^{Q_1}_{\mathbf{n_p},\mathbf{n_q}}+H^{Q_1}_{\mathbf{n_p},\mathbf{n_q}}+T^{Q_2}_{\mathbf{n_p},\mathbf{n_q}}+H^{Q_2}_{\mathbf{n_p},\mathbf{n_q}})\mathbf{Q}_{vec}=\mathbf{0},
\end{equation}
where $T^{Q_l}_{\mathbf{n_p},\mathbf{n_q}}$ and $H^{Q_l}_{\mathbf{n_p},\mathbf{n_q}},$ for $l=1,2$ are the block Toeplitz and block Hankel matrices, respectively defined as $$T^{Q_1}_{\mathbf{n_p},\mathbf{n_q}}(r,s) = (T_{r-s}),~~T^{Q_2}_{\mathbf{n_p},\mathbf{n_q}}(r,s) = (H_{r-s}),~~ H^{Q_1}_{\mathbf{n_p},\mathbf{n_q}}(r,s) = (T_{r+s}),~~ H^{Q_2}_{\mathbf{n_p},\mathbf{n_q}}(r,s) = (H_{r+s}),$$ 
for $r=n_{p_x}+1,n_{p_x}+2,\ldots, n_{p_x}+n_{q_x}+1, s=0,1,2,\ldots,n_{q_x}$.

Observe that, each block matrix in \eqref{2DComputingQ} has $(n_{q_x}+1)\times(n_{q_x}+1)$ blocks of $(n_{q_y}+1)\times(n_{q_y}+1)$ size. Therefore, after removing the last row of the elements sitting in the last row of the block matrix system \eqref{2DComputingQ}, we get a homogeneous system with $(n_{q_x}+1)\times(n_{q_y}+1)-1$ equations in  $(n_{q_x}+1)\times(n_{q_y}+1)$ unknowns which always has a non-trivial solution.

A unique bivariate Pad\'e-Chebyshev approximant can be obtained if the coefficient matrix given in \eqref{2DComputingQ} with Chebyshev series coefficients is of full rank and in the case of rank deficiency one can choose any vector of the basis of the kernel of the matrix as a solution $\textbf{Q}_{vec}$. Similar to the univariate case, we use MATLAB $null$ subroutine to compute the basis of the kernel of the matrix and we have chosen the last column of the kernel space as $\textbf{Q}_{vec}$.

On expanding the index $i$ and $j$, we can write \eqref{2Dinhomosystem} in the matrix form as follows
\begin{equation}\label{2DComputingP}
\mathbf{P}_{vec}=\dfrac{1}{4}\left(\begin{bmatrix}
A_0 \\
A_1 & A_0\\
\vdots \\
A_{n_{q_x}} & A_{n_{q_x}-1}  & \cdots & A_0\\
\vdots \\
A_{n_{p_x}} & A_{n_{p_x}-1}  & \cdots & A_{n_{p_x}-n_{q_x}}
\end{bmatrix}+\begin{bmatrix}
A_0 & A_1  & \cdots & A_{n_{q_x}}\\
A_1 & A_2  & \cdots & A_{n_{q_x}+1}\\
\vdots \\
A_{n_{q_x}} & A_{n_{q_x}+1}  & \cdots & A_{2n_{q_x}}\\
\vdots \\
A_{n_{p_x}} & A_{n_{p_x}+1} & \cdots & A_{n_{p_x}+n_{q_x}}
\end{bmatrix}+\begin{bmatrix}
0 & 0  & \cdots & 0 \\
0 & A_0  & \cdots & A_{n_{q_x}-1} \\
\vdots \\
0 & 0 &  \cdots & A_0 \\
\vdots \\
0 & 0  & \cdots & 0
\end{bmatrix}\right)\mathbf{Q}_{vec},
\end{equation}
where each element $A_m$ is the sum of two Toeplitz and a Hankel matrices of size $(n_{p_y}+1)\times (n_{q_y}+1)$ given by 
\begin{equation}\label{eq:1Dinhomoshortmatrix}
A_m=\begin{bmatrix}
c_{m,0} \\
c_{m,1} & c_{m,0}\\
\vdots \\
c_{m,n_{q_y}} & c_{m,n_{q_y}-1}  & \cdots & c_{m,0}\\
\vdots \\
c_{m,n_{p_y}} & c_{m,n_{p_y}-1}  & \cdots & c_{m,n_{p_y}-n_{q_y}}
\end{bmatrix}+\begin{bmatrix}
c_{m,0} & c_{m,1}  & \cdots & c_{m,n_{q_y}}\\
c_{m,1} & c_{m,2}  & \cdots & c_{m,n_{q_y}+1}\\
\vdots \\
c_{m,n_{q_y}} & c_{m,n_{q_y}+1}  & \cdots & c_{m,2n_{q_y}}\\
\vdots \\
c_{m,n_{p_y}} & c_{m,n_{p_y}+1} & \cdots & c_{m,n_{p_y}+n_{q_y}}
\end{bmatrix}+\begin{bmatrix}
0 & 0  & \cdots & 0 \\
0 & c_{m,0}  & \cdots & c_{m,n_{q_y}-1} \\
\vdots \\
0 & 0 &  \cdots & c_{m,0} \\
\vdots \\
0 & 0  & \cdots & 0
\end{bmatrix}.
\end{equation}
for $m = 0,1,2,\ldots,n_{p_x}+n_{q_x}$.

In our numerical simulations, we use approximated Chebyshev series coefficients \eqref{eq:bivarChebcoeffapprx} computed by using the two-dimensional Gauss-Chebyshev quadrature formula to compute bivariate Pad\'e-Chebyshev approximation and denote it by $\mathsf{R}^{n_x,n_y,M}_{\mathbf{n_p},\mathbf{n_q}}(x,y)$. Here $n_x$ and $n_y$ are the number of Chebyshev points in $x$ and $y$-direction, respectively, we use in the quadrature formula. 
\subsection{Piecewise 2D Pad\'{e}-Chebyshev Approximation}\label{sec:Pi2DPC}
Let $f$ be a bivariate bounded function defined on a rectangular domain  $[a_x,b_x]\times[a_y,b_y]$. For given two-tuples $\mathbf{n}=(n_x,n_y)\ge\mathbf{1}$, $\mathbf{n_p} = (n_{p_x},n_{p_y})\ge\mathbf{n_q} = (n_{q_x},n_{q_y})\ge\mathbf{1}$ and  $\mathbf{N}=(N_x,N_y)\ge\mathbf{1}$, where $n_x,n_y,n_{p_x},n_{q_x},n_{p_y},n_{q_y},N_x,N_y$ are integers, construct a piecewise bivariate Pad\'e-Chebyshev approximation of $f$ as follows:
\begin{enumerate}
\item Discretize the domain into $N_x$ cells in $x$-direction and $N_y$ cells in $y$-direction, denoted by $I_{j_x} = [a_{j_x},b_{j_x}],j_x = 0,1, \ldots, N_x-1$ and $I_{j_y} = [a_{j_y},b_{j_y}],j_y = 0,1,\ldots, N_y-1$, respectively, as explained in the step 1 of Section \ref{sec:Pi2DC}.
\item For $j_x = 0,1,2,\ldots,N_x-1$ and $j_y = 0,1,2,\ldots,N_y-1$, obtain the truncated Chebyshev series in each sub rectangle $I_{j_x}\times I_{j_y}$, as explained in the Step 3 and 4 of Section \ref{sec:Pi2DC} with $d_x^{j_x} = n_{p_x}+2n_{q_x}$ and $d_y^{j_y}=n_{p_y}+2n_{q_y}$.
\item  Define the piecewise bivariate Pad\'e-Chebyshev approximation of $f$ in the rectangular domain $I_x\times I_y$, where $I_x= [a_x,b_x]$ and $I_y = [a_y,b_y]$, with respect to the given partition as
\begin{equation}
\mathsf{R}_{\mathbf{n_p},\mathbf{n_q}}^{\mathbf{n},\mathsf{M},N_x,N_y}(x,y) = 
	\begin{cases}
	\mathsf{R}_{n_p,n_q,n}^{\mathsf{M},I_{0_x}\times I_{0_y}}(x,y) &\quad\text{if } (x,y) \in [a_{0_x},b_{0_x})\times [a_{0_y},b_{0_y})\\
	\mathsf{R}_{n_p,n_q,n}^{\mathsf{M},I_{0_x}\times I_{1_y}}(x,y) &\quad\text{if } (x,y) \in [a_{0_x},b_{0_x})\times [a_{1_y},b_{1_y})\\
	\vdots \\
	\mathsf{R}_{n_p,n_q,n}^{\mathsf{M},I_{0_x}\times I_{N_y}}(x,y) &\quad\text{if } (x,y) \in [a_{0_x},b_{0_x})\times [a_{N_y-1},b_{N_y-1})\\
	\vdots \\
	\mathsf{R}_{n_p,n_q,n}^{\mathsf{M},I_{N_x-1}\times I_{0_y}}(x,y) &\quad\text{if } (x,y) \in [a_{N_x-1},b_{N_x-1})\times [a_{0_y},b_{0_y})\\
	\mathsf{R}_{n_p,n_q,n}^{\mathsf{M},I_{N_x-1}\times I_{1_y}}(x,y) &\quad\text{if } (x,y) \in [a_{N_x-1},b_{N_x-1})\times [a_{1_y},b_{1_y})\\
	\vdots \\
	\mathsf{R}_{n_p,n_q,n}^{\mathsf{M},I_{N_x-1}\times I_{N_y}}(x,y) &\quad\text{if } (x,y) \in [a_{N_x-1},b_{N_x-1})\times [a_{N_y-1},b_{N_y-1}),
	\end{cases}
	\end{equation}
where $\mathsf{R}_{n_p,n_q,n}^{\mathsf{M},I_{j_x}\times I_{j_y}}(x,y)$, for $(x,y)\in I_{j_x}\times I_{j_y}$, denotes the Maehly based Pad\'e-Chebyshev approximant of $f|I_{j_x}\times I_{j_y} , j_x = 0, 1,\ldots, N_x-1$ and $j_y = 0, 1,\ldots, N_y-1$.
\end{enumerate}
Hence $\mathsf{R}_{\mathbf{n_p},\mathbf{n_q}}^{\mathbf{n},\mathsf{M},N_x,N_y}(x,y)$ is the desired piecewise bivariate Pad\'e-Chebyshev approximant of the function $f$ on the domain $I_x\times I_y$ which we referred to as Pi2DPC.

\section{Numerical Comparison}\label{sec:2Dnumericalresults}
\begin{figure}[H]
\centering
\includegraphics[height=8cm,width=8cm]{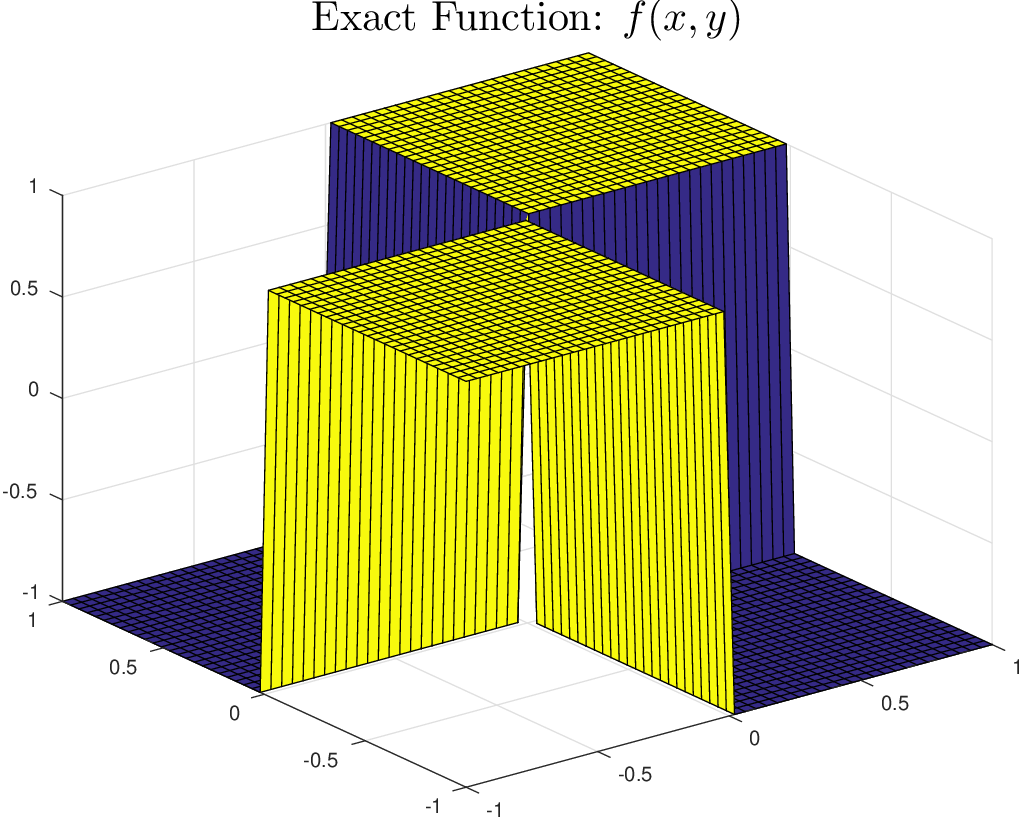}
\caption{Depicts the function $f(x,y)= sign(4xy)$ on the unit square $D=[-1,1]^2$.}
\label{fig:2Dsign}
\end{figure}

\begin{figure}[H]
\centering
\includegraphics[height=7cm,width=7cm]{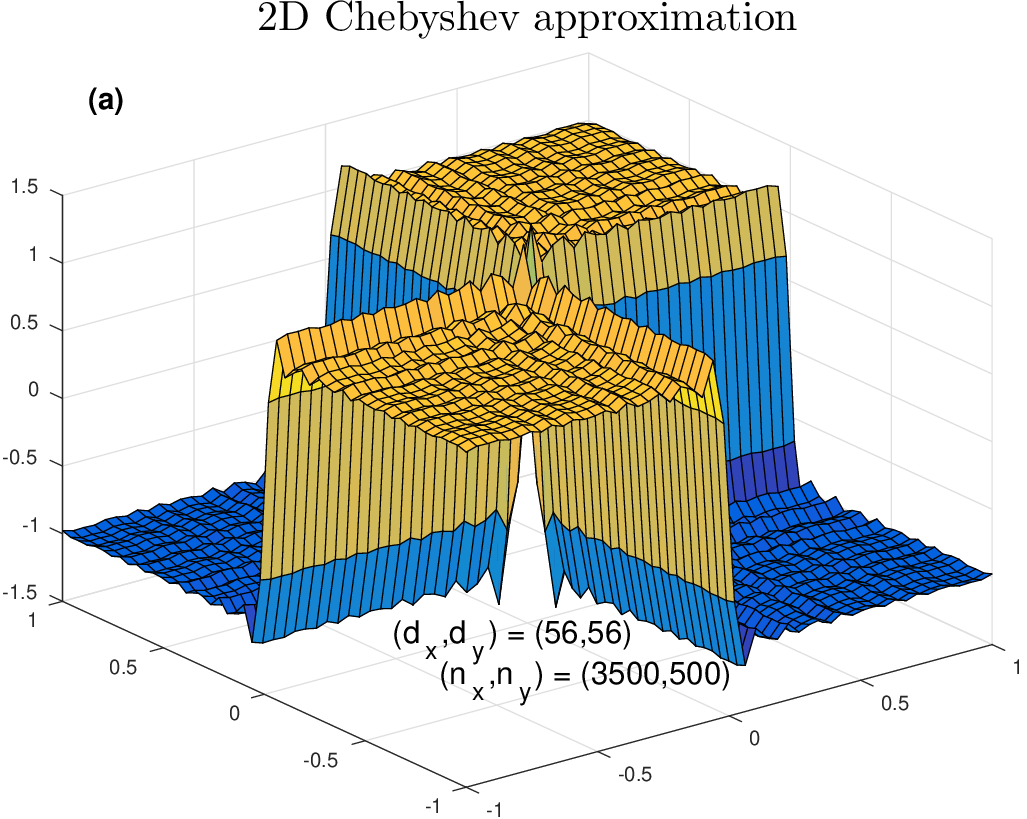}
\includegraphics[height=7cm,width=7cm]{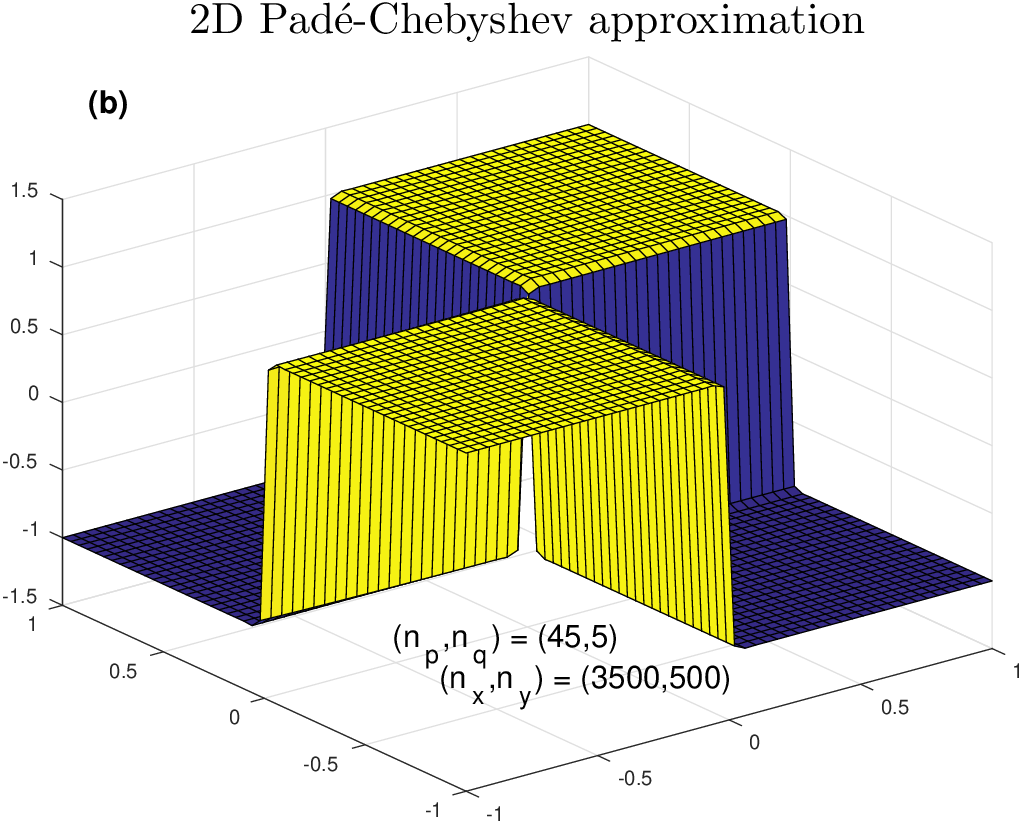}\\
\includegraphics[height=7cm,width=7cm]{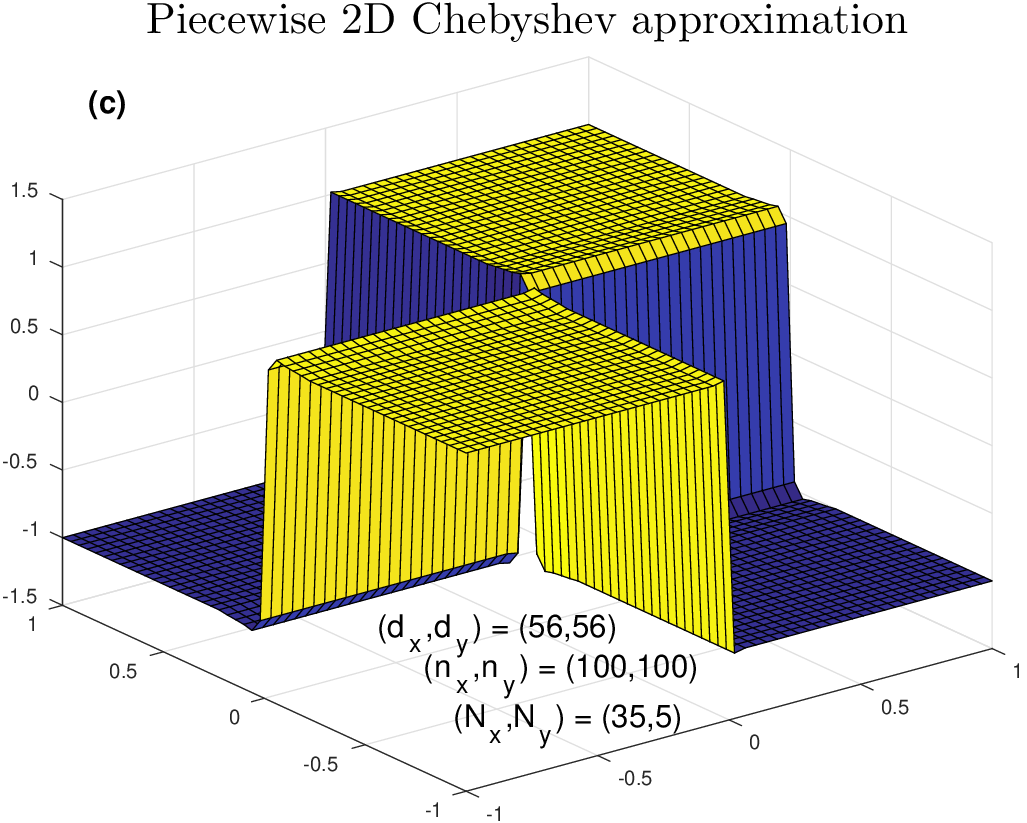}
\includegraphics[height=7cm,width=7cm]{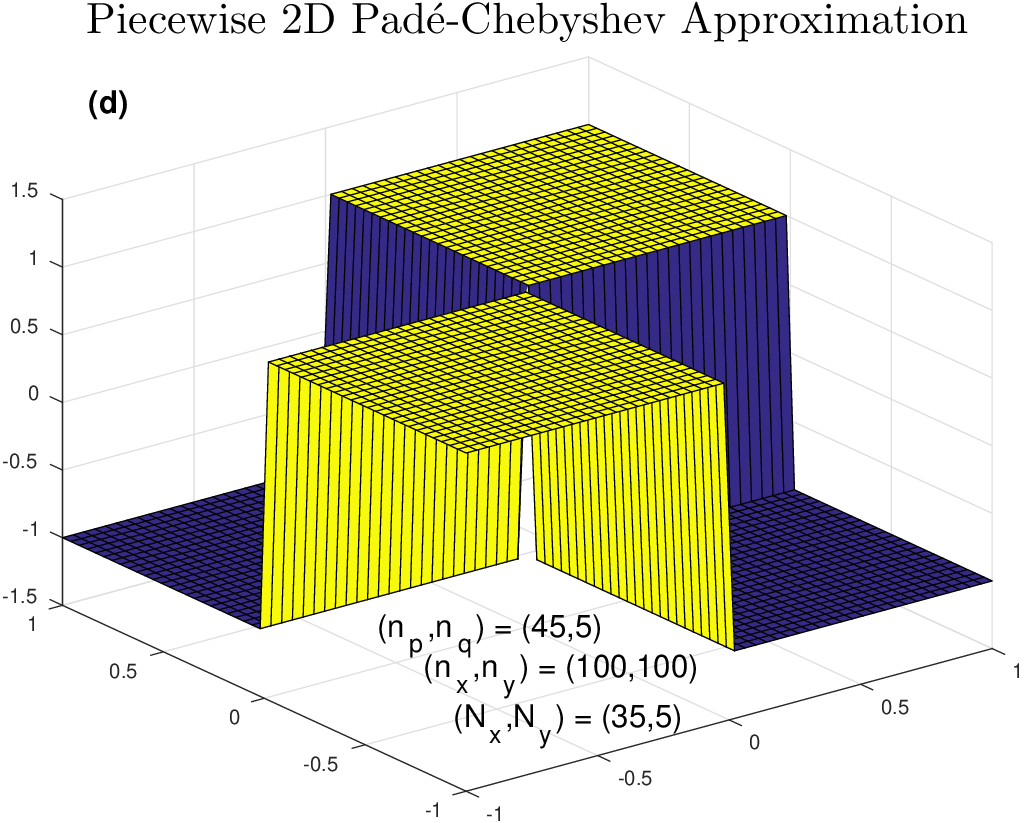}
\caption{Depicts \textbf{(a)} global bivariate Chebyshev, \textbf{(b)} global bivariate Pad\'e-Chebyshev, \textbf{(c)} Pi2DC and \textbf{(d)} Pi2DPC approximants of the function $f(x,y)= sign(4xy), (x,y)\in D$.}
\label{fig:2DCPCsgn}
\end{figure}

In this section, we study the performance of the proposed Pi2DPC technique while approximating bivariate piecewise smooth functions. We consider few test examples to validate the method and present numerical evidence that the proposed algorithm captures singularity curves in two-dimensional functions without a visible Gibbs phenomenon.
\begin{example}\label{ex:2Dsgn}
	Let us consider the following real valued piecewise constant function  
	\begin{equation}\label{eq:2Dsgn}
	z=f(x,y) = sign(4xy), ~~~~~~ (x,y) \in D.
	\end{equation}
The function is symmetric with respect to $x$ and $y$ axes and involves jump discontinuity along the straight lines at $x=0$ and $y=0$, as shown in Figure \ref{fig:2Dsign}. We approximate the function $f$ using bivariate Chebyshev (see Figure \ref{fig:2DCPCsgn}\textbf{(a)}), global bivariate Pad\'e-Chebyshev (see Figure \ref{fig:2DCPCsgn}\textbf{(b)}), Pi2DC (see Figure \ref{fig:2DCPCsgn}\textbf{(c)}) and Pi2DPC (see Figure \ref{fig:2DCPCsgn}\textbf{(d)}) approximation methods. 
\end{example}

The proposed Pi2DC and Pi2DPC algorithms, developed in Section \ref{sec:Pi2DC} and \ref{sec:Pi2DPC}, respectively, are performed to approximate $f$ by fixing $\mathbf{N}=(35,5)$ and we use $\mathbf{n}=(100,100)$ data points in Gauss quadrature formula to approximate the Chebyshev series coefficients in each sub domain $I_{j_x}\times I_{j_y},j_x = 0,1, \ldots, N_x-1,j_y = 0,1,\ldots, N_y-1$. We choose $\mathbf{n_p}=(45,45)$, $\mathbf{n_q}=(5,5)$ and $(d_x,d_y)=(56,56)$ in each sub domain to perform the piecewise algorithms. To perform the global algorithms, defined in Section \ref{sec:2Dchebyshevapprox} and \ref{sec:2DMPCapprox}, we fix the parameters $\mathbf{n_p}=(45,45)$, $\mathbf{n_q}=(5,5)$ and $(d_x,d_y)=(56,56)$.  We use $\mathbf{n}=(35\times100,5\times100)$ points in Gauss quadrature formula to approximate the Chebyshev series coefficients.

The comparison demonstrated in Figure \ref{fig:2DCPCsgn} shows that the four approximants are well in agreement with the exact function $f$, in the smooth segments. Moreover, the Pi2DPC approximant captures the singularity curves sharply (see Figure \ref{fig:2DCPCsgn}\textbf{(d)}) with negligible Gibbs oscillations in the vicinity. It is evident from the Figure \ref{fig:2DCPCsgn}\textbf{(a)}, \textbf{(b)} and \textbf{(c)}  that the global approximants and the Pi2DC exhibits Gibbs phenomenon and are able to capture the singularity curves accurately with less resolution. However, we can clearly see the importance of a rational approximation in approximating bivariate non-smooth functions. 

\begin{figure}[H]
\centering
\includegraphics[height=7cm,width=7cm]{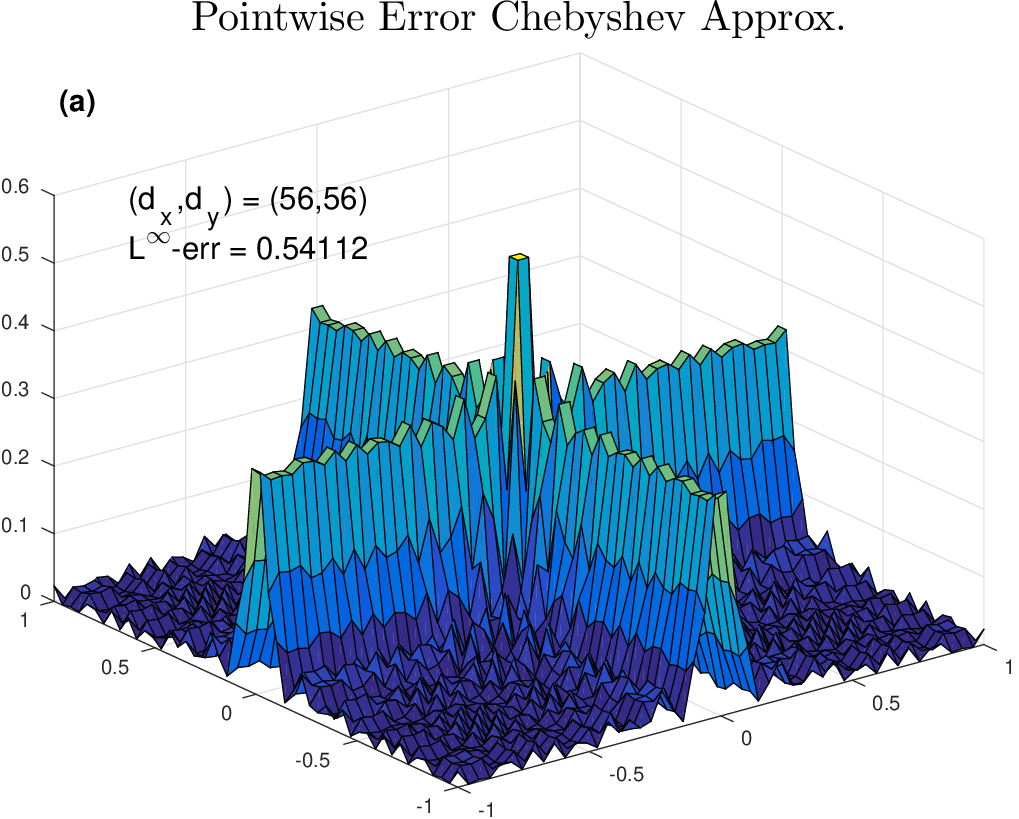}
\includegraphics[height=7cm,width=7cm]{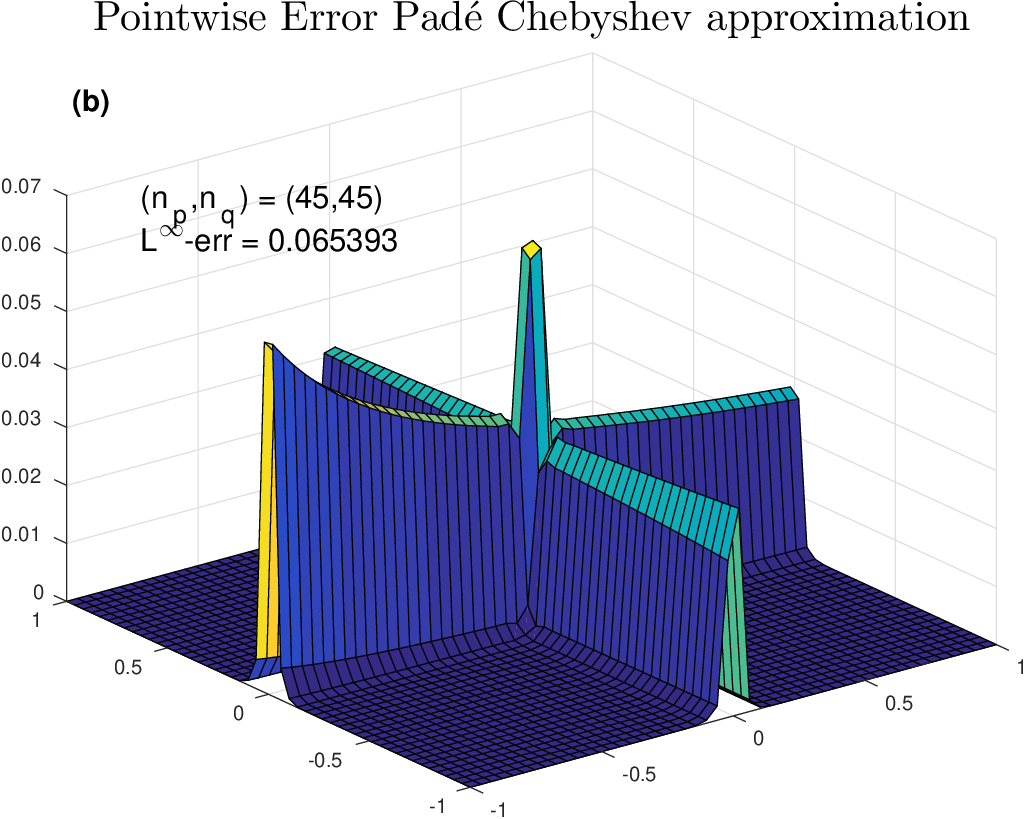}\\
\includegraphics[height=7cm,width=7cm]{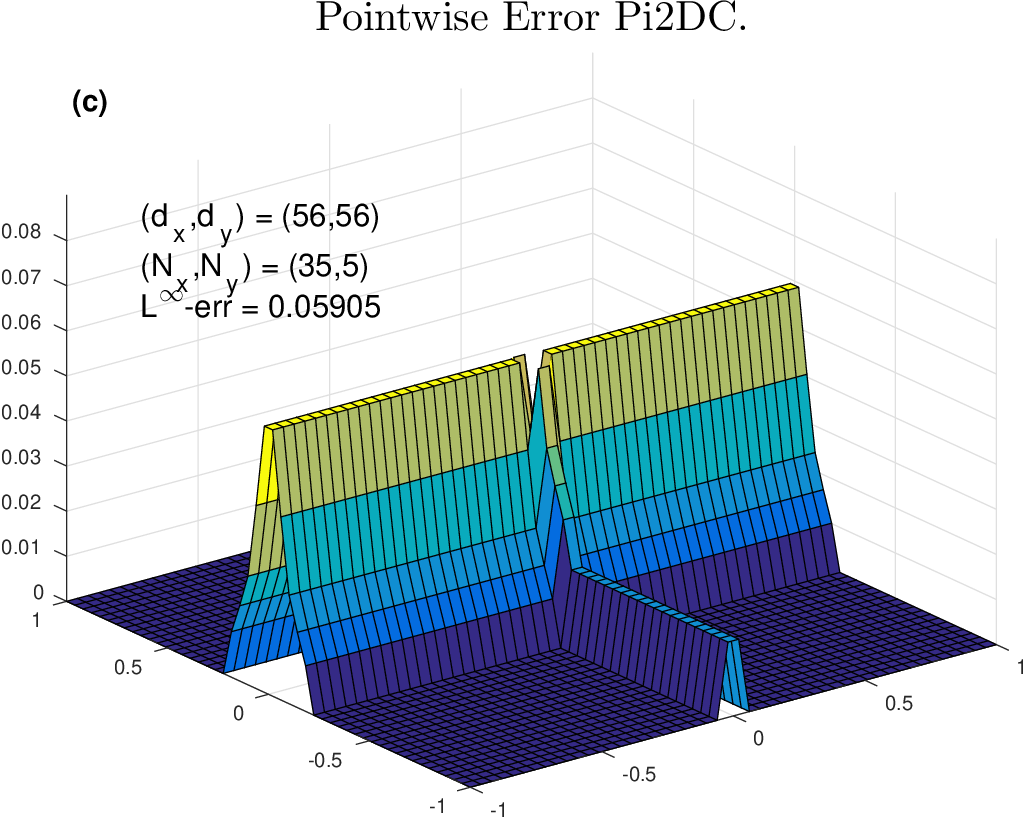}
\includegraphics[height=7cm,width=7cm]{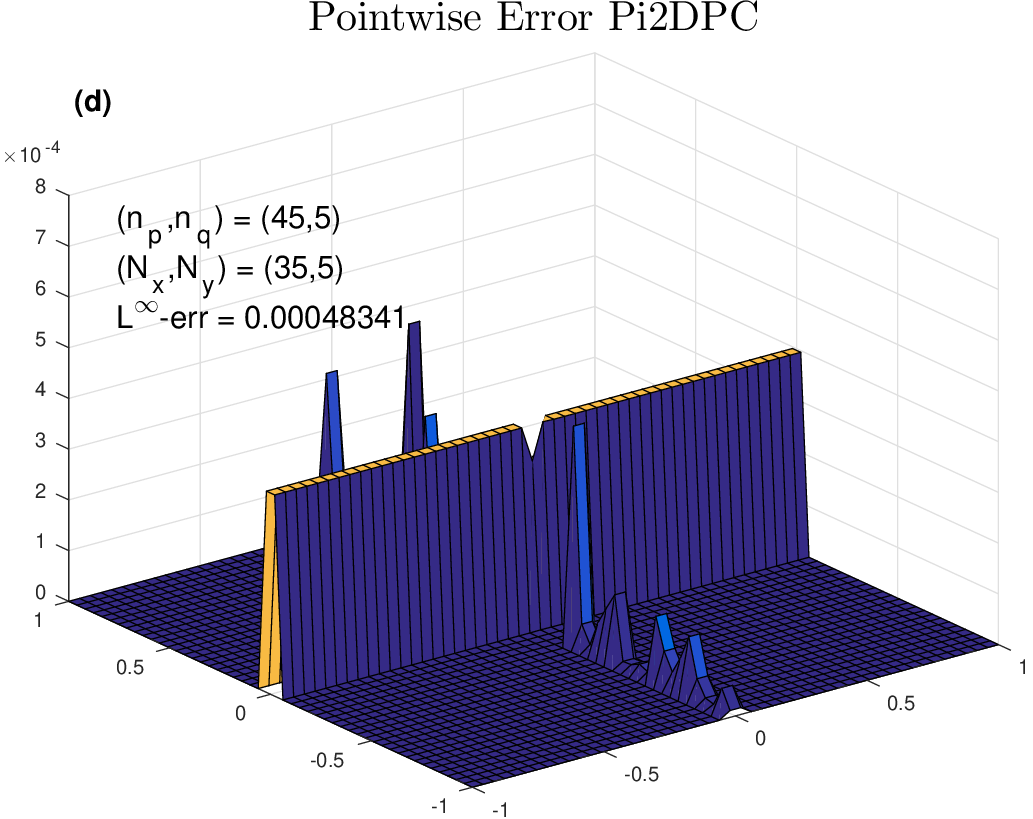}
\caption{The corresponding pointwise error plots are depicted with $L^\infty$-errors.}
\label{fig:2DCPCsgnerr}
\end{figure}

Since, often in applications exact Chebyshev coefficients are not known, we use quadrature formula to approximate these coefficients for numerical experiments. Computation of approximated Chebyshev coefficients for a given function can be done efficiently. The most expensive step of Pi2DPC algorithm is the inversion of a matrix of size $(n_{q_x}+1)^2\times(n_{q_y}+1)^2$ in each sub rectangle. However, we observed that to capture the singularity curves accurately in a bivariate piecewise smooth function using Pi2DPC method one need to choose a small value for $n_{q_x}$ and $n_{q_y}$. 
The pointwise error graphs of the approximants is depicted in Figure \ref{fig:2DCPCsgnerr} alongside the corresponding $L^\infty$-errors. The maximum errors mentioned in Figure \ref{fig:2DCPCsgnerr}\textbf{(a)} and \ref{fig:2DCPCsgnerr}\textbf{(b)} indicates that the global Pad\'e-Chebyshev approximant is $10$ times better than the global Chebyshev series approximant in the vicinity of singularity. A better approximation property of rational functions makes Pad\'e-Chebyshev approximation method a right choice for approximating bivariate piecewise smooth functions. Figure \ref{fig:2DCPCsgnerr}\textbf{(c)} and \ref{fig:2DCPCsgnerr}\textbf{(d)} shows a clear significance of the piecewise implementation of these approximants as the $L^\infty$-error of Pi2DC approximant is $10$ times better than its global counterpart. It is evident from Figure \ref{fig:2DCPCsgnerr}\textbf{(d)} that the Pi2DPC approximant is the best representation of the function $f$ among four. Additionally, the $L^\infty$-error is the least for Pi2DPC approximant.
  
\begin{figure}[H]
\includegraphics[height=6cm,width=14cm]{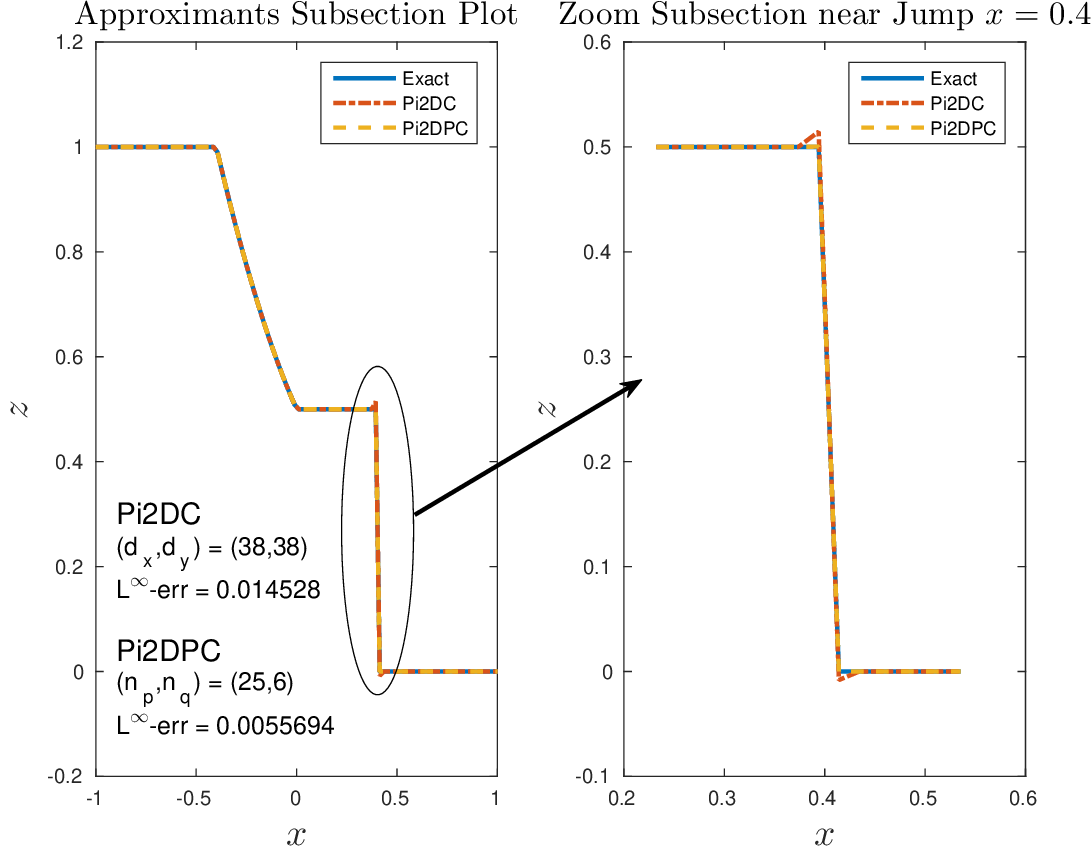}
\caption{Comparision of the slice of the surface of Pi2DC and Pi2DPC approximants of $h(x,y)$ along with the zoomed plot.}
	\label{fig:2Dcrosssecdisc2}
\end{figure}

\begin{example}\label{ex:disc2}
Let us consider the following function as our second test case which has a jump discontinuity and a low order singularity (where $h(x,y)$ is continuous but $h_x(x,y)$ is discontinuous along a curve) along a straight line parallel to $y$-axis.
	\begin{equation}\label{eq:disc2}
	z=h(x,y) = 
	\begin{cases}
	1 &\quad\text{if } x \in [-1,-0.4),y\in [-1,1]\\
	x^2-\frac{17}{20}x+\frac{1}{2} &\quad\text{if } x \in [-0.4,0),y\in [-1,1]\\
	\dfrac{1}{2} &\quad\text{if } x \in [0,0.4),y\in [-1,1]\\
	0 &\quad\text{if } t \in [0.4,1],y\in [-1,1],
	\end{cases}
	\end{equation}
	We perform the proposed Pi2DC and Pi2DPC techniques to approximate the above defined function. The cross-section of this function along $x$-axis (see Figure \ref{fig:2Dcrosssecdisc2}) look-alike the solution of \textit{Sod shock tube} problem in one space dimension. The solution of such problems develop shock (jump) in finite time and the main challenge in approximating these solutions is the development of the Gibbs oscillations.
\end{example}
\begin{figure}[H]
	\centering
	\includegraphics[height=7cm,width=5.2cm]{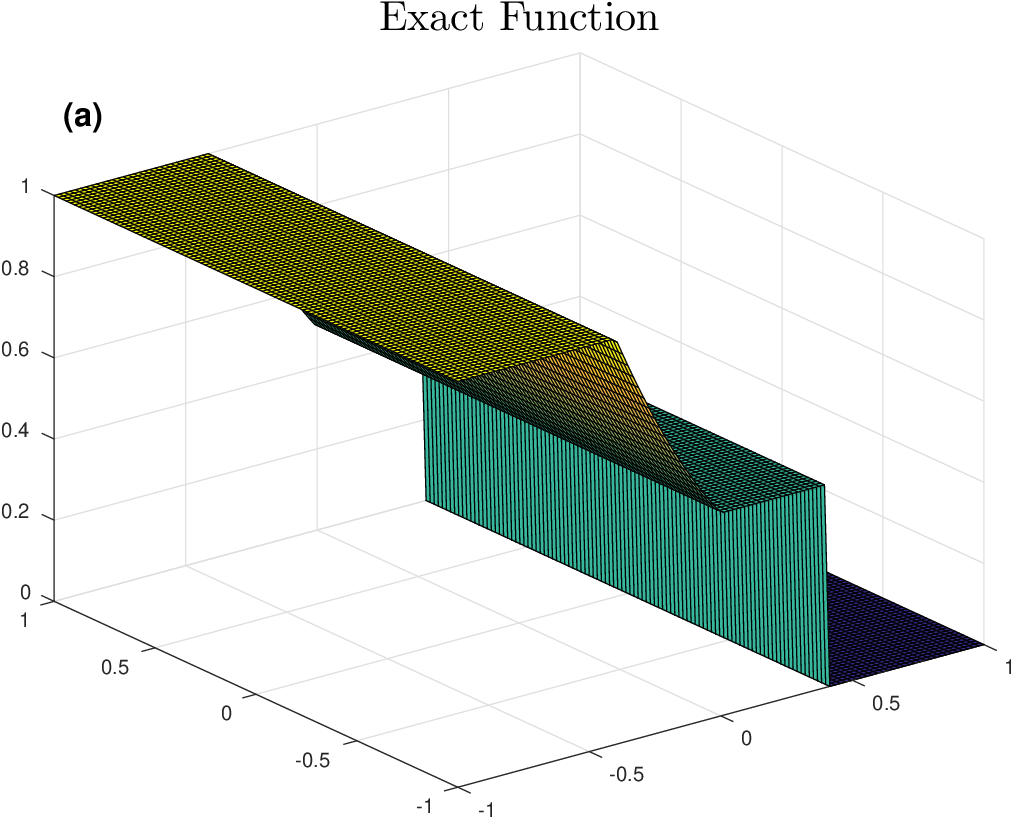}
	\includegraphics[height=7cm,width=5.2cm]{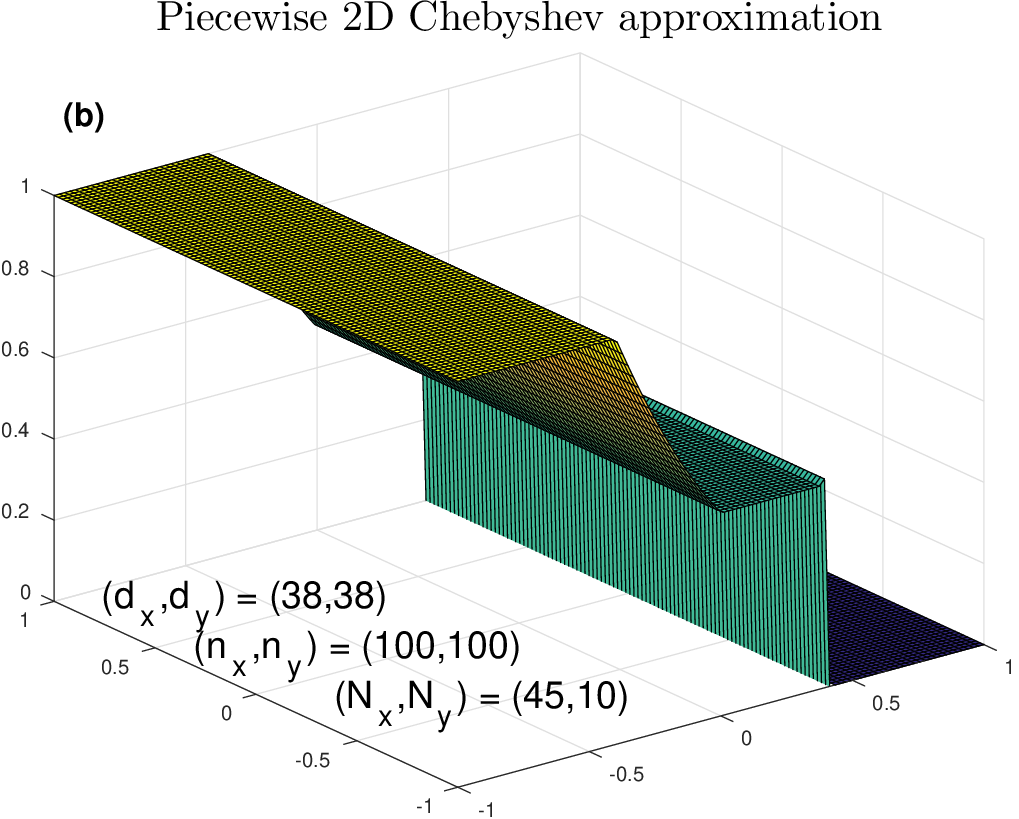}
	\includegraphics[height=7cm,width=5.2cm]{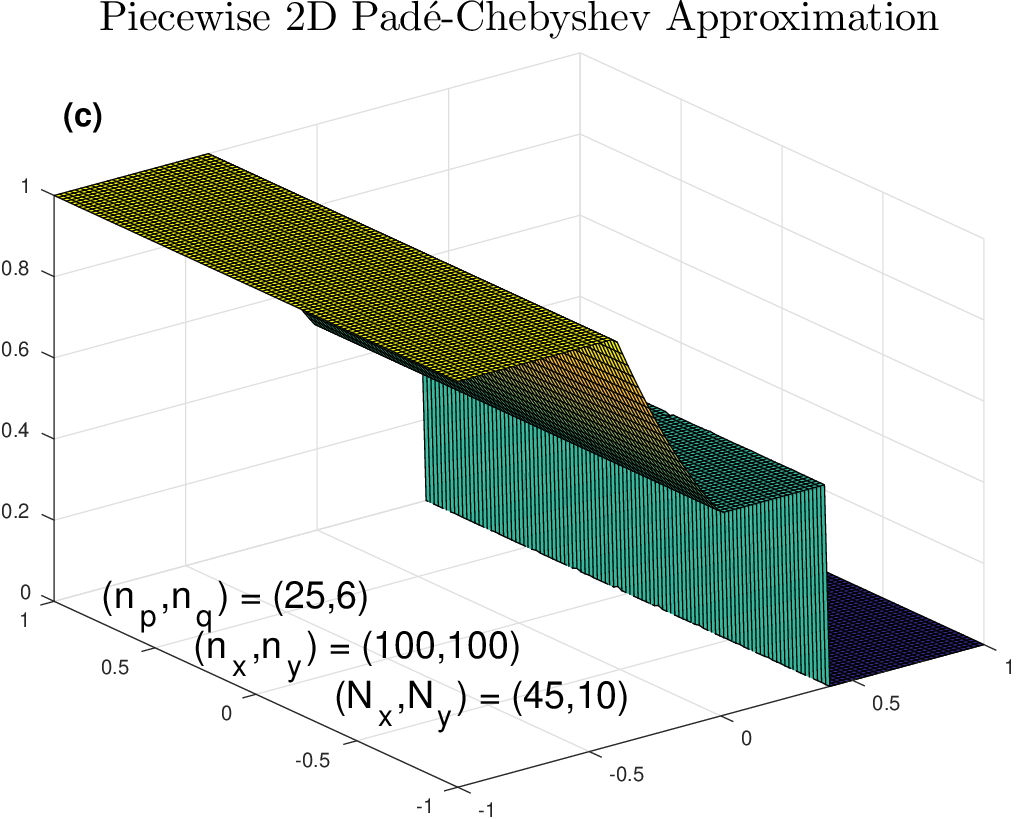}
	\caption{\textbf{(a)} Depicts the surface plot of the function defined in \eqref{eq:disc2} on the unit square $[-1,1]^2$. \textbf{(b)} depicts Pi2DC approximant, and \textbf{(c)} depicts Pi2DPC approximant of $h(x,y)$.}
	\label{fig:2Ddisc2}
\end{figure}
The function $h(x,y)$ involves a lower order singularity at $x=-0.4$ and a jump discontinuity at $x=0.4$ in the $z$-plane as shown in Figure \ref{fig:2Ddisc2}\textbf{(a)}. To perform the piecewise algorithms, we fixed the parameter $\mathbf{N}=(45,10)$, $\mathbf{n}=(100,100)$, $\mathbf{n_p}=(25,25)$, $\mathbf{n_q}=(6,6)$ and $(d_x,d_y)=(38,38)$ in each sub domain $I_{j_x}\times I_{j_y}$ for $j_x = 0,1, \ldots, N_x-1$, and $j_y = 0,1,\ldots, N_y-1$. The comparison demonstrated in Figure \ref{fig:2Ddisc2} shows that both the approximants are well in agreement with the exact function $f$ in the smooth segments. However as expected, the Pi2DPC approximant is able to capture the jump (see Figure \ref{fig:2Dcrosssecdisc2}) with negligible Gibbs oscillations in the vicinity of the singularities. To achieve this accuracy, Pi2DPC algorithm solves a block linear system of size $49\times 48$ for $450$ times. We use MATLAB $null(A)$  subroutine to solve the homogeneous linear system. The Pi2DPC algorithm generates Figure \ref{fig:2Ddisc2}\textbf{(c)} in less than 4 seconds and Pi2DC algorithm takes around 2 seconds to generate Figure \ref{fig:2Ddisc2}\textbf{(b)} on a standard desktop. 


Figure \ref{fig:2Dcrosssecdisc2} depicts a comparison between the approximants Pi2DC and Pi2DPC in a cross-section along $x$-axis. This profile of the surface \eqref{eq:disc2} look similar to the Sod shock tube problem in one space dimensional at some non-zero time. From these plots we can see that both Pi2DC and Pi2DPC approximants captures the point singularity at $x=-0.4$ with high resolution. However, as expected from previous example and one dimensional case, we can clearly see a lack in accuracy of the Pi2DC approximant at the jump $x=0.4$ whereas Pi2DPC approximant captures the jump accurately with almost no oscillations. Moreover, the $L^\infty$-error in Pi2DPC approximant listed in Figure \ref{fig:2Dcrosssecdisc2} indicates the significance of a rational approximation.


\section{Conclusion}
In this article, we proposed a piecewise Pad\'e-Chebyshev(PiPC) algorithm for approximating non smooth functions on a bounded interval. We have tested the same against functions possessing a jump discontinuity as well as a point(lower order) singularity, to demonstrate numerically, a faster rate of convergence (as compared to the global Pad\'e-Chebyshev approximation) of the approximants in the vicinity of singularities as the number of cells and the degrees of the numerator and denominator polynomials approach infinity. The PiPC algorithm does not require a prior knowledge of the location and type of singularities present in the original function and is well suited for implementation using a non uniform mesh, thereby offering greater flexibility. Also, the proposed piecewise algorithm works well for functions which involves singularities at irrational points.  Extending the PiPC to two dimensions, we then propose the piecewise bivariate Pad\'e-Chebyshev algorithm (Pi2DPC) to approximate bivariate piecewise smooth functions. Pi2DPC too, like its univariate counterpart PiPC, does not rely on any prior information about the location and type of singularities present in the target function. Numerical experiments carried out with functions possessing different order of singularities along certain curves clearly demonstrate the effectiveness of the Pi2DPC algorithm in capturing the singularity curves sharply without any visible Gibbs phenomenon. Finally, the results obtained using Pi2DPC are compared with those resulting from piecewise and global bivariate Chebyshev approximation.
\bibliography{bibfile_thesis}

\begin{thebibliography}{10}

\bibitem{aka_bas-19a}
S.~Akansha and S.~Baskar.
\newblock Adaptive pad\'e-chebyshev type approximation to piecewise smooth
  functions, 2019.

\bibitem{dmi_yos-12a}
Dmitry B. and Yosef Y.
\newblock Algebraic fourier reconstruction of piecewise smooth functions.
\newblock {\em Mathematics of Computation}, 81(277):277--318, 2012.

\bibitem{bak-geo_65a}
G.~A. Baker, Jr.
\newblock The theory and application of the {P}ad\'e approximant method.
\newblock In {\em Advances in {T}heoretical {P}hysics, {V}ol. 1}, pages 1--58.
  Academic Press, New York, 1965.

\bibitem{bak-geo-gra-pet_96a}
G.~A. Baker, Jr. and G.-M. Peter.
\newblock {\em Pad\'e approximants}, volume~59 of {\em Encyclopedia of
  Mathematics and its Applications}.
\newblock Cambridge University Press, Cambridge, second edition, 1996.

\bibitem{bas_73a}
N.~K. Basu.
\newblock On double chebyshev series approximation.
\newblock {\em SIAM Journal on Numerical Analysis}, 10(3):496--505, 1973.

\bibitem{bat-15a}
D.~Batenkov.
\newblock Complete algebraic reconstruction of piecewise-smooth functions from
  {F}ourier data.
\newblock {\em Math. Comp.}, 84(295):2329--2350, 2015.

\bibitem{ber-12a}
S.N. Bernstein.
\newblock {Sur l’ordre de la meilleure approximation des fonctions continues
  par des polynômes de degré donné}.
\newblock {\em Mem. Acad. R. Belg}, 13:1--104, 1912.

\bibitem{ber-14a}
S.N. Bernstein.
\newblock {Sur la meilleure approximation de |x| par des polynômes de degrés
  donnés}.
\newblock {\em Acta Mathematica}, 37:1--57, 1914.

\bibitem{boy-01a}
J.~P. Boyd.
\newblock {\em Chebyshev and {F}ourier spectral methods}.
\newblock Dover Publications, Inc., Mineola, NY, second edition, 2001.

\bibitem{che_66a}
E.~W. Cheney.
\newblock {\em Introduction to approximation theory}.
\newblock McGraw-Hill Book Co., New York-Toronto, Ont.-London, 1966.

\bibitem{che_98a}
E.~W. Cheney.
\newblock {\em Introduction to approximation theory}.
\newblock AMS Chelsea Publishing, Providence, RI, 1998.
\newblock Reprint of the second (1982) edition.

\bibitem{dev_lor-93a}
R.~A. DeVore and G.~G. Lorentz.
\newblock {\em Constructive approximation}, volume 303 of {\em Grundlehren der
  Mathematischen Wissenschaften [Fundamental Principles of Mathematical
  Sciences]}.
\newblock Springer-Verlag, Berlin, 1993.

\bibitem{dri_for-01a}
T.~A. Driscoll and B.~Fornberg.
\newblock A {P}ad\'e-based algorithm for overcoming the {G}ibbs phenomenon.
\newblock {\em Numer. Algorithms}, 26(1):77--92, 2001.

\bibitem{eck-93a}
K.~S. Eckhoff.
\newblock Accurate and efficient reconstruction of discontinuous functions from
  truncated series expansions.
\newblock {\em Math. Comp.}, 61(204):745--763, 1993.

\bibitem{gee_ban-97a}
J.~Geer and N.~S. Banerjee.
\newblock Exponentially accurate approximations to piece-wise smooth periodic
  functions.
\newblock {\em Journal of Scientific Computing}, 12(3):253--287, Sep 1997.

\bibitem{hes_kab_lur-06a}
J.~S. Hesthaven, S.~M. Kaber, and L.~Lurati.
\newblock Pad\'e-{L}egendre interpolants for {G}ibbs reconstruction.
\newblock {\em J. Sci. Comput.}, 28(2-3):337--359, 2006.

\bibitem{kab_mad-05a}
S.~M. Kaber and Y.~Maday.
\newblock Analysis of some {P}ad\'e-{C}hebyshev approximants.
\newblock {\em SIAM J. Numer. Anal.}, 43(1):437--454, 2005.

\bibitem{lan-16a}
C.~Lanczos.
\newblock {\em Discourse on Fourier Series}.
\newblock Classics in Applied Mathematics. Society for Industrial and Applied
  Mathematics, 2016.

\bibitem{lev-20a}
D.~Levin.
\newblock Reconstruction of piecewise-smooth multivariate functions from
  fourier data, 2020.

\bibitem{lit-03a}
G.~L. Litvinov.
\newblock Error autocorrection in rational approximation and interval
  estimates. [{A} survey of results].
\newblock {\em Cent. Eur. J. Math.}, 1(1):36--60, 2003.

\bibitem{lut_74a}
C.~H. Lutterodt.
\newblock A two-dimensional analogue of {P}ad\'{e} approximant theory.
\newblock {\em J. Phys. A}, 7:1027--1037, 1974.

\bibitem{lut_75a}
C.~H. Lutterodt.
\newblock Addendum to: ``{A} two-dimensional analogue of {P}ad\'{e} approximant
  theory'' ({J}. {P}hys. {A} {\bf 7} (1974), 1027--1037).
\newblock {\em J. Phys. A}, 8:427--428, 1975.

\bibitem{mae_60a}
H.~J. Maehly.
\newblock Rational approximations for transcendental functions.
\newblock In {\em Information processing}, pages 57--62. UNESCO, Paris; R.
  Oldenbourg, Munich; Butterworths, London, 1960.

\bibitem{mar_pau_wei_hen_las-21a}
S.~Marx, E.~Pauwels, T.~Weisser, D.~Henrion, and J.~B. Lasserre.
\newblock {Semi-algebraic Approximation Using Christoffel–Darboux Kernel}.
\newblock {\em Constructive Approximation}, 2021.

\bibitem{mas-han-03a}
J.~C. Mason and D.~C. Handscomb.
\newblock {\em Chebyshev polynomials}.
\newblock Chapman \& Hall/CRC, Boca Raton, FL, 2003.

\bibitem{min_kab_don-07a}
M.~S. Min, S.~M. Kaber, and W.~S. Don.
\newblock Fourier-{P}ad\'e approximations and filtering for spectral
  simulations of an incompressible {B}oussinesq convection problem.
\newblock {\em Math. Comp.}, 76(259):1275--1290, 2007.

\bibitem{new-64a}
D.~J. Newman.
\newblock {Rational approximation to $| x| $.}
\newblock {\em Michigan Mathematical Journal}, 11(1):11 -- 14, 1964.

\bibitem{nur_89a}
G.~N{\"u}rnberger.
\newblock {\em Approximation by spline functions}.
\newblock Springer-Verlag, Berlin, 1989.

\bibitem{pac_pla_tre-10a}
R.~Pach\'{o}n, R.~B. Platte, and L.~N. Trefethen.
\newblock Piecewise-smooth chebfuns.
\newblock {\em IMA J. Numer. Anal.}, 30(4):898--916, 2010.

\bibitem{phi_03a}
G.~M. Phillips.
\newblock {\em Interpolation and approximation by polynomials}.
\newblock CMS Books in Mathematics/Ouvrages de Math\'ematiques de la SMC, 14.
  Springer-Verlag, New York, 2003.

\bibitem{riv-74a}
T.~J. Rivlin.
\newblock {\em The {C}hebyshev polynomials}.
\newblock Wiley-Interscience [John Wiley \& Sons], New York-London-Sydney,
  1974.
\newblock Pure and Applied Mathematics.

\bibitem{run-85a}
C.~Runge.
\newblock {Zur Theorie der Eindeutigen Analytischen Functionen}.
\newblock {\em Acta Mathematica}, pages 229--244, 1885.

\bibitem{sho_73a}
I.~J. Schoenberg.
\newblock {\em Cardinal spline interpolation}.
\newblock Society for Industrial and Applied Mathematics, Philadelphia, Pa.,
  1973.
\newblock Conference Board of the Mathematical Sciences Regional Conference
  Series in Applied Mathematics, No. 12.

\bibitem{tam_lop_hes-12a}
A.~L. Tampos, J.~E.~C. Lope, and J.~S. Hesthaven.
\newblock Accurate reconstruction of discontinuous functions using the singular
  {P}ad\'e-{C}hebyshev method.
\newblock {\em IAENG Int. J. Appl. Math.}, 42(4):242--249, 2012.

\bibitem{you-chu_09a}
T.~Youtian and Chuanqing Gu.
\newblock A two-dimensional matrix {P}ad\'{e}-type approximation in the inner
  product space.
\newblock {\em J. Comput. Appl. Math.}, 231(2):680--695, 2009.

\end{thebibliography}
\end{document}